\documentclass[12pt]{article}
\usepackage{amssymb}
\usepackage{amsmath}
\numberwithin{equation}{section}

\newcommand{\bea}{\begin{eqnarray}}
\newcommand{\eea}{\end{eqnarray}}
\newcommand{\beano}{\begin{eqnarray*}}
\newcommand{\eeano}{\end{eqnarray*}}
\newcommand{\nonu}{\nonumber \\}

\newcommand{\hs}[1]{\hspace{#1 mm}}
\newcommand{\eps}{\epsilon}

\newcommand{\lda}{\lambda}
\newcommand{\Lda}{\Lambda}

\newcommand{\bvarphi}{\overline{\varphi}}
\newcommand{\Phidag}{\Phi^\dagger}
\newcommand{\phidag}{{\phi^\dagger}}

\newcommand{\bL}{{\overline{L}}}
\newcommand{\bct}{{\overline{\ct}}}
\newcommand{\lambdadag}{\lambda^\dagger}
\newcommand{\gvph}{{\mbox{\boldmath{$\varphi$}}}}

\newcommand{\gch}{{\mbox{\boldmath{$\chi$}}}}
\newcommand{\gps}{{\mbox{\boldmath{$\psi$}}}}
\newcommand{\gdlt}{{\mbox{\boldmath{$\delta$}}}}
\newcommand{\gld}{{\mbox{\boldmath{$\lambda$}}}}

\newcommand{\gf}{\mathbf f}
\newcommand{\bfg}{\mathbf g}
\newcommand{\gF}{\mathbf F}
\newcommand{\gG}{\mathbf G}

\newcommand{\gA}{\mathbf A}
\newcommand{\edag}{e^\dagger}
\newcommand{\Adag}{A^\dagger}
\newcommand{\gAdag}{\gA^\dagger}
\newcommand{\tAdag}{\tilde{A}^\dagger}
\newcommand{\tgAdag}{\tilde{\gA}^\dagger}
\newcommand{\adag}{a^\dagger}
\newcommand{\bino}[2]{\left(
          \begin{array}{c}#1\\#2\end{array}\right)}

\newcommand{\ca}{\mbox{$\cal{A}$}}

\newcommand{\cd}{\mbox{$\cal{D}$}}

\newcommand{\cf}{\mbox{${\cal F}$}}
\newcommand{{\cg}}{\mbox{$\cal{G}$}}
\newcommand{\ch}{\mbox{$\cal{H}$}}

\newcommand{\cl}{\mbox{$\cal{L}$}}

\newcommand{\calr}{\mbox{$\cal{R}$}}
\newcommand{\ct}{\mbox{$\cal{T}$}}

\newcommand{\prt}{\partial}

\newcommand{\mb}[1]{\hs{4}\mbox{#1}\hs{4}}
\newcommand{\half}{\frac{1}{2}}

\newtheorem{theo}{Theorem}[section]
\newtheorem{prop}[theo]{Property}

\newtheorem{defi}[theo]{Definition}

\newtheorem{lem}[theo]{Lemma}

\newcommand{\prf}{\underline{Proof:}\ }
\newcommand{\finprf}{\null \hfill {\rule{5pt}{5pt}}\\ \null}
\newcommand{\ie}{{\it i.e.}\ }
\newcommand{\CC}{\mbox{${\mathbb C}$}}
\newcommand{\RR}{\mbox{${\mathbb R}$}}
\newcommand{\ZZ}{\mbox{${\mathbb Z}$}}
\newcommand{\1}{\mbox{\hspace{.0em}1\hspace{-.24em}I}}
\newcommand{\II}{\mbox{${\mathbb I}$}}

\newcommand{\pp}{{\mathbf p}}
\newcommand{\qq}{{\mathbf q}}

\newcommand{\NP}[1]{Nucl.\ Phys.\ {\bf #1}}

\newcommand{\CMP}[1]{Commun.\ Math.\ Phys.\ {\bf #1}}
\newcommand{\JMP}[1]{J.\ Math.\ Phys.\ {\bf #1}}

\newcommand{\IJMP}[1]{Int. J.\ Mod.\ Phys.\ {\bf #1}}
\newcommand{\PR}[1]{Phys.\ Rev.\ {\bf #1}}

\newcommand{\lcom}{{[\hspace{-.3ex}[}}
\newcommand{\rcom}{{]\hspace{-.3ex}]}}

\begin{document}
\renewcommand{\thefootnote}{\fnsymbol{footnote}}
\newpage
\pagestyle{empty} \setcounter{page}{0}
%%%%%%%%%%%%%%%%%%%%%%%%%%%%%%%%
%%%%%  HEADINGS POUR DRAFT  %%%%%%

\markright{\today\dotfill DRAFT\dotfill } 
%\pagestyle{myheadings}

%%%%%%%%%%%%%%%%%%%%%%%%%%%%%%%%

%%%%%%%%%%%%%%%%%%%%%%%%%%%%%%%%%%%%%%%%%%%%%%%%%
%%%%%%%%%%%%%%%%%% LOGO LAPTH - DEBUT  %%%%%%%%%%%%%%%
%%%%%%%%%%%%%%%%%%%%%%%%%%%%%%%%%%%%%%%%%%%%%%%%%
\newcommand{\LAP}{LAPTH}
\def\logo{{\bf {\huge LAPTH}}}

\centerline{\logo}

\vspace {.3cm}

\centerline{{\bf{\it\Large Laboratoire d'Annecy-le-Vieux de
Physique Th\'eorique}}}

\centerline{\rule{12cm}{.42mm}}
%%%%%%%%%%%%%%%%%%%%%%%%%%%%%%%%%%%%%%%%%%%%%%%%%%
%%%%%%%%%%%%%%%%% LOGO LAPTH  - FIN %%%%%%%%%%%%%%
%%%%%%%%%%%%%%%%%%%%%%%%%%%%%%%%%%%%%%%%%%%%%%%%%%

\vspace{20mm}

\begin{center}

{\LARGE  {\sffamily Lax pair and super-Yangian symmetry of the \\[.21cm]
non-linear super-Schr\"odinger equation}}

\vspace{10mm}

{\large V. Caudrelier \footnote{caudreli@lapp.in2p3.fr} and E.
Ragoucy\footnote{ragoucy@lapp.in2p3.fr}\\[.21cm]
 Laboratoire de Physique Th{\'e}orique \LAP\footnote{UMR 5108
  du CNRS, associ{\'e}e {\`a} l'Universit{\'e} de
Savoie.}\\[.242cm]
 LAPP, BP 110, F-74941  Annecy-le-Vieux Cedex, France. }
\end{center}
\vfill\vfill

\begin{abstract}
We consider a version of the non-linear Schr\"odinger equation with
$M$ bosons and $N$ fermions. We first solve the classical and quantum
versions of this equation, using a super-Zamolodchikov-Faddeev (ZF)
algebra. Then we prove that the hierarchy associated to this model
admits a super-Yangian $Y(gl(M|N))$ symmetry. We exhibit the
corresponding (classical and quantum) Lax pairs. Finally, we
construct explicitly the super-Yangian generators, in terms of the
canonical fields on the one hand, and in terms of the ZF algebra
generators on the other hand. The latter construction uses the
well-bred operators introduced recently.
\end{abstract}

\vfill
\rightline{\tt math.QA/0306115}
\rightline{\LAP-984/03}
\rightline{June 2003}

\newpage
%\pagestyle{plain}
%%%%%%%%%%%%%%%%%%%%%%%%%%%%%%%%
%%%%%  FIN PAGE DE TITRE  %%%%%%
%%%%%%%%%%%%%%%%%%%%%%%%%%%%%%%%
\setcounter{footnote}{0}

\section{Introduction}

The nonlinear Schr{\"o}dinger (NLS) equation is one of the most
studied systems in quantum integrable systems (for a review, see
e.g. \cite{Gut}), and its simplest (scalar) version played an
important r\^ole in the development of the (Quantum) Inverse
Scattering Method \cite{Zhak}. As usual in quantum integrable
systems, its integrability relies on the existence of an infinite
dimensional symmetry algebra. In integrable systems, natural
candidates for such algebras are the celebrated quantum groups
associated to (affine) Lie algebras, or the Yangians. Indeed, it
is known~\cite{MuWa} that the quantum NLS model with spin
$\frac{1}{2}$ fermions and repulsive interaction on the line has a
Yangian symmetry $Y(sl(2))$. More generally, its vectorial
version, based on $N$-component bosons or on $N$-component
fermions, was shown to possess a $Y(gl(N))$ symmetry \cite{MRSZ}.
The integrability can also be grounded on the existence of an
infinite series of mutually commuting Hamiltonians, which thus
generates a whole hierarchy of equations. In the case of scalar
NLS equation, the hierarchy contains well-known models, such as
the modified KdV equation.

It was natural to seek a supersymmetric version (including both
bosons and fermions) of these models which admits
 the super-Yangian based on superalgebras $gl(M|N)$ as symmetry
algebra. Different versions of such a generalization were already
proposed, from the simple boson-fermion systems related to NLS
\cite{sup1,sup2}, or superfields formulation \cite{superchp} of
NLS, up to more algebraic studies of these models
\cite{KP-NLS,superDS}. The difficulty with such generalizations is
to keep the fundamental notion of integrability while allowing for
the existence of supersymmetry. Even when some of the suggested
supersymmetric systems were shown to pass some integrability
conditions \cite{sup3}, the status of such models remained not
clearly established, and one is still looking for e.g. their Lax
presentation or their underlying infinite dimensional symmetry
algebra.

Another $\ZZ_{2}$-graded version of NLS was introduced by Kulish
\cite{Z2Ku}, the fields being super-matrix valued and thus
associated to both fermions and bosons. However, only the finite
interval was studied, using the Thermodynamical Bethe Ansatz (see
also \cite{FPZ}), and the explicit quantum solutions are not
known. The symmetry (super)algebra is also lacking in this
presentation.

The aim of this article is to present a "super-vectorial" version
(close to the matricial version introduced by Kulish) of the NLS
model on the infinite line which includes $M$ bosons \textbf{and}
$N$ fermions fields. The advantage of this version relies on its
manifest integrability and the existence of quantum canonical
solutions, that we will explicitly construct using a super-ZF
algebra
 (section \ref{sec2}). Indeed,
these solutions can be associated to a whole hierarchy of mutually
commuting Hamiltonians, as it should be for an integrable model.
It also admits, as we will show (section \ref{sec3}), a Lax presentation
both at classical and quantum level
(without using a superfield formalism). As usual, the Lax pair
presentation allows to recover the hierarchy of our super-NLS equation.
Finally, this super NLS
hierarchy possesses a super-Yangian symmetry and we will construct it,
 both using the quantum canonical solutions or the super-ZF generators
 (section \ref{sec4}).

\section{Non-linear super-Schr\"odinger equation\label{sec2}}

\subsection{The usual Non-linear Schr\"odinger equation}
The NLS equation reads
\begin{equation}
\label{CVNLS}
\left(i\partial_t+\partial^2_x\right)\phi_i(x,t)=2g\phidag^j(x,t)\phi_j(x,t)\phi_i(x,t),~i=1,...,N
\mb{with}g>0
\end{equation}
where summation over repeated indices is understood. It is obtained from the (time-independent) Hamiltonian
\begin{equation}
\label{classical-hamiltonian}
H(\phi_i,\phidag_j)=\int_{-\infty}^{\infty}dx\,
\left(\partial_x\phidag^j(x)\partial_x\phi_j(x)+g\phidag^i(x)\phidag^j(x)\phi_j(x)\phi_i(x)\right)
\end{equation}
using the Hamiltonian equation of motion $\partial_t F=\{H,F\}$,
valid for any functional
$F(\phi_i,\phidag_j)$,
where the Poisson Bracket (PB) is canonically associated to $\phi$
and $\phidag$.

A solution \textit{\`{a} la
Rosales} \cite{Ros} can be written as follows
\begin{equation}
\label{solution-rosales1}
\phi_i(x,t)=\sum_{n=0}^{\infty}(-g)^n\phi_i^{(n)}(x,t)\,,\qquad
g>0
\end{equation}
with
$$
\phi_i^{(n)}(x,t)=\int_{\RR^{2n+1}}d^n\pp d^{n+1}\qq\,
\lambda^{k_1}(p_1)\cdots\lambda^{k_n}(p_n)\lambda_{k_n}(q_n)\cdots\lambda_{k_1}(q_1)
\lambda_{i}(q_0)\frac{e^{i\Omega_{n}(x,t;\pp,\qq)}}{Q_{n}(\pp,\qq,0)}
$$
\begin{eqnarray}
\Omega_{n}(x,t;\pp,\qq) &=& \sum\limits_{j=0}^n(q_j x-q^2_j t)-
\sum\limits_{i=1}^n(p_i x-p^2_i t)\nonu
Q_{n}(\pp,\qq,\varepsilon) &=&
\prod\limits_{i=1}^n(p_i-q_{i-1}+i\varepsilon)(p_i-q_i+i\varepsilon)\nonu
d^n\pp d^{n+1}\qq &=& \prod_{i=1\atop j=0}^n
\frac{dp_i}{2\pi}\frac{dq_j}{2\pi}
\label{solution-rosales2}
\end{eqnarray}
where we have denoted $\pp=(p_{1},\ldots,p_{n})$,
$\qq=(q_{0},\ldots,q_{n})$.

The Rosales solution is fundamental since its structure is
preserved upon quantization \cite{Dav} and we shall see below that
this result survives when one includes fermions. The NLS equation
and its hierarchy admit the Yangian $Y(gl(N))$ as symmetry, the
explicit construction of its generators was given in \cite{MuWa}
(for $sl(2)$, in terms of canonical fields) and \cite{MRSZ} (for
$sl(N)$, in terms of the ZF generators). A Lax pair formulation
can be found in \cite{skly79, Thac} (for NLS equation) and in
\cite{PZ,PWZ} (for its vectorial generalisation).

\subsection{Classical Non-linear super-Schr\"odinger equation}
We consider a generalized version of the NLS eq. which includes both
bosons and fermions. Due to the use of auxiliary spaces (see
appendix), the corresponding eq. will formally look like the original
one, but let us insist that the present version is a "supersymmetric"
version of it. While the similarities allow us to build the solution
of the Non-linear super-Schr\"odinger equation, the differences will
appear for instance in the nature of the symmetry algebra (see below).

We define
$\Phi(x)=\sum_{j=1}^{M+N}\phi_{j}(x)e_j$,
where $e_{j}$ is a $(M+N)$-column vector in the auxiliary space (see
appendix)
and summation is understood for repeated indices.
$\phi_j,~j=1,...,M$ and $\phi_j,~j=M+1,...,M+N$ are the bosonic and
fermionic components respectively. For convenience, we set
$K=M+N$.
We shall also need adjoints of the fields
\begin{equation}
\Phidag(x)=\phidag_{i}(x)\edag_i~,~~~x\in\RR
\end{equation}
The Hamiltonian reads:
\begin{equation}
H(\Phi,\Phidag)=\int_{-\infty}^{\infty}dx\,
\left(\,\partial_x\Phidag(x)\partial_x\Phi(x)+
g\,\left(|\Phi(x)|^{2}\right)^2\,\right)
\end{equation}
or in components:
\begin{equation}
H(\Phi,\Phidag)=\int_{-\infty}^{\infty}dx\,
\left(\partial_x\phidag^{j}(x)\partial_x\phi_j(x)+
g\,\phidag^{j}(x)\phidag^{k}(x)\phi_k(x)\phi_j(x)\right)
\label{defi-hamiltonian}
\end{equation}
The canonical Poisson brackets for the
basic fields $\Phi(x)$, $\Phidag(y)$ with corresponding components
$\phi_i(x)$, $\phidag_j(y)$ take the following form
\begin{eqnarray}
   \label{global-bracket}
   \{\Phi_1(x),\Phidag_2(y)\} &=& i\delta_{12}\delta(x-y)
   =-\{\Phidag_2(y),\Phi_1(x)\}
   \quad\mb{(globally)}\quad\\
    \label{comp-bracket}
    \{\phi_j(x),\phidag_k(y)\} &=& i\delta_{jk}\delta(x-y)
    =-(-1)^{[j][k]} \{\phidag_k(y),\phi_j(x)\}
   \mb{(in components)}\qquad
\end{eqnarray}
The field $\Phi(x,t)$ of components $\phi_i(x,t)$
satisfies the following Hamiltonian equation of motion which we
call the classical Nonlinear super-Schr\"odinger (NLSS)
equation.
\begin{eqnarray}
\label{CNLSS-eq1} i\partial_t\Phi(x,t) &=&\!
-\partial^2_x\Phi(x,t)+2g|\Phi(x,t)|^2\,
\Phi(x,t)\qquad\qquad\quad\mb{(globally)}\quad\\
\label{CNLSS-eq2} i\partial_t\phi_j(x,t) &=&\!
-\partial^2_x\phi_j(x,t)+2
g\,(\phidag_{k}(x,t)\phi_k(x,t))\,\phi_j(x,t) \quad\mb{(in
components)}\qquad
\end{eqnarray}
 These equations are simply derived from the Hamiltonian
equations of motion $\partial_t\Phi(x,t)=\{H,\Phi(x,t)\}$ and
$\partial_t\phi_i(x,t)=\{H,\phi_i(x,t)\}$. The equations of motion
are (formally) the same as the usual ones and the solution
\textit{\`{a} la Rosales} (\ref{solution-rosales1}),
(\ref{solution-rosales2}) is still valid in our case:

\begin{theo}
The solution of the classical NLSS equation (\ref{CNLSS-eq2}) is given by
\begin{eqnarray}
\phi_j(x,t)&=&\sum_{n=0}^{\infty}(-g)^n\phi^{(n)}_j(x,t)\
\mb{where}\\
\phi^{(n)}_j(x,t)&=&\int_{\RR^{2n+1}}d^n\pp
d^{n+1}\qq \,\sum_{k_{1},\ldots,k_{n}=1}^{K}
\lambdadag_{k_{1}}(p_1)\cdots\lambdadag_{k_{n}}(p_n)\lambda_{k_{n}}(q_n)\cdots
\lambda_{k_{1}}(q_1)\lambda_j(q_0)\nonu &
&\times\frac{e^{i\Omega_{n}(x,t;\pp,\qq)}}{Q_{n}(\pp,\qq,0)}
\end{eqnarray}
using the same notations as in (\ref{solution-rosales2}).
\end{theo}
\prf Substituting into NLSS equation, it amounts to the following
identity being satisfied
$$\sum_{j=0}^n q_j^2-\sum_{i=1}^n p_i^2-\left(\sum_{j=0}^n q_j-\sum_{i=1}^n p_i\right)^2
=-2\sum_{c=1}^{n-1}\sum_{a=1}^c(p_{a+1}-q_a)(p_{c+1}-q_{c+1})$$
which is readily seen to hold.\finprf Note that, due to the
$\ZZ_{2}$-graded tensor product, the ordering of the
$\lda^\dagger$'s and of the $\lda$'s respectively matters.

\subsection{Quantizing NLSS\label{sect-quant}}

\subsubsection{Graded ZF algebra}
We write a graded version of the ZF algebra \cite{ZZ}-\cite{F},
using auxiliary spaces and entities containing  bosonic
and fermionic components (see appendix):
\begin{equation}
\gA(k)    =a_{i}(k)e_i
~~~~\text{ and }~~~~ \gAdag(k)
        =a^{\dagger}_{i}(k)\edag_i,\quad k\in\RR
\end{equation}

\begin{defi}
The graded ZF algebra reads
\begin{eqnarray}
\label{ZF-algebra1}
\gA_1(k_1)\gA_2(k_2)&=&R_{21}(k_2-k_1)\gA_2(k_2)\gA_1(k_1)\\
\label{ZF-algebra2}
\gAdag_1(k_1)\gAdag_2(k_2)&=&\gAdag_2(k_2)\gAdag_1(k_1)R_{21}(k_2-k_1)\\
\label{ZF-algebra3}
\gA_1(k_1)\gAdag_2(k_2)&=&\gAdag_2(k_2)R_{12}(k_1-k_2)\gA_1(k_1)+
\gdlt_{12}\delta(k_1-k_2)
\end{eqnarray}
where
\begin{equation}
R_{12}(k)=\frac{k\1\otimes\1-igP_{12}}{k+ig}\label{R-matrix}
\end{equation}
is the R-matrix for
the super-Yangian $Y(gl(M|N))\equiv Y(M|N)$, and $P_{12}$ is the
super-permutation operator:
\begin{equation}
P_{12}=\sum_{i,j=1}^{K}(-1)^{[j]}E_{ij}\otimes E_{ji} \label{defP12}
\end{equation}
\end{defi}
Note that for even vectors $u$, $v$ and even matrices $B$, $C$ (as
defined in appendix), one has $P_{12}\, (u\otimes v)
= v\otimes u$ and $P_{12}\, (B\otimes C)\, P_{12}= C\otimes B$.

The $R$-matrix has the following useful properties
\begin{eqnarray}
&&R_{21}(k)=R_{12}(k)\\
&&R_{12}(k_1-k_2)R_{21}(k_2-k_1)=\1\otimes\1\\
&&R^\dagger_{12}(k_1-k_2)=R_{21}(k_2-k_1)
\end{eqnarray}

{F}or quantities of definite $\ZZ_2$-grade, we define
 their {super-commutator} by
\begin{equation}
[\hspace{-2pt}[B,C]\hspace{-2pt}]=BC-(-1)^{[B][C]}CB
\end{equation}
Then, after some calculations, one shows that the component version of the ZF algebra reads ($j,k=1,\ldots,K$):
\begin{eqnarray}
\hspace{-2.4ex}
\left[\hspace{-4pt}\left[\rule{0ex}{2.4ex}a_j(k_1),a_k(k_2)\right]\hspace{-4pt}\right]&=&
\frac{-ig}{k_2-k_1+ig}\,\left(\rule{0ex}{2.4ex}a_j(k_2)a_k(k_1)+(-1)^{[j][k]}a_k(k_2)a_j(k_1)\right)
\quad\\
\left[\hspace{-4pt}\left[\adag_j(k_1),\adag_k(k_2)\right]\hspace{-4pt}\right]&=&
\frac{-ig}{k_2-k_1+ig}\,\left(\adag_j(k_2)\adag_k(k_1)+(-1)^{[j][k]}\adag_k(k_2)\adag_j(k_1)\right)\\
\left[\hspace{-4pt}\left[a_j(k_1),\adag_k(k_2)\right]\hspace{-4pt}\right]&=&
\frac{-ig}{k_1-k_2+ig}\left((-1)^{[j][k]}\adag_k(k_2)a_j(k_1)+\delta_{jk}\,\sum_{\ell
=1}^{K} \adag_{\ell}(k_2)a_{\ell}(k_1)\right)\nonu
&&+\delta_{jk}\,\delta(k_1-k_2)
\end{eqnarray}
Note that these relations ensure the existence of a PBW basis,
generated by the monomials having $\adag$'s on the left of the
$a$'s,
 the $a$'s on one hand, and the $\adag$'s on the other hand
being ordered according to the
magnitude of the "impulsions" $k_{j}$.
\subsubsection{Fock representation}

The previous algebra can be represented on a Fock space, which is
most useful for our quantization of NLSS, and we follow the basic
ideas of \cite{Dav} (further developed in e.g. \cite{LM} and
\cite{GLM}). A detailed presentation of the graded version when
$M=N=1$ has been given in \cite{CR}. The general case follows the
same lines, so that we just sketch the results, referring to
\cite{CR} for more details about the $\ZZ_{2}$-graded case.

We introduce
 $\cf_R=\bigoplus_{n=0}^{\infty}\ch_R^n$ where
 $\ch_R^0=\CC$,
$$
\ch_R^1=\left\{\gvph(p)=\sum_{j=1}^{K}\varphi_{j}(p)e_j
\mbox{ s.t. }\varphi_j\in L^2(\RR),~j=1,\ldots,{K}\right\}
\equiv KL^2(\RR)
$$
and for $n\ge
2$:
\begin{eqnarray*}
\hspace{-2.4ex} \ch_R^n &=& \Big\{ \gvph_{1...
n}(p_1,...,p_n)=\sum_{i_1,...,i_n=1}^{K}
\varphi_{i_1,...,i_n}(p_1,...,p_n)(e_{i_1}\otimes\cdots\otimes
e_{i_n})
\\
&&\mb{s.t.}\varphi_{i_1,...,i_n}\in
L^2(\RR^n),~i_1,...,i_n=1,\ldots,{K} \mb{and}\\
\lefteqn{\hspace{-4.2em}\gvph_{1... i,i+1...
n}(p_1,...,p_i,p_{i+1},...,p_n)=R_{i,i+1}(p_i-p_{i+1}) \gvph_{1...
i+1,i... n}(p_1,...,p_{i+1},p_i,...,p_n)\Big\}}
\end{eqnarray*}

 There exists a (vacuum) vector $\Omega\in\cd$ which is
cyclic with respect to $\gAdag(k)$ and annihilated by $\gA(k)$.

The scalar product which we define below on $\ch_R^n$
provides the usual $L^2$ topology and $\cf_R$ is the completed
vector space over $\CC$ for this topology.

The sesquilinear form $\langle~,~\rangle$ defined on
$\ch_R^n\times\ch_R^n,~n\ge 1$ by
\begin{eqnarray}
\label{scalar-product} \langle\gvph,\gps\rangle
&=&\int_{\RR^n}d^np~
\gvph^{\dagger}_{1...n}(p_1,...,p_n)\gps_{1...n}(p_1,...,p_n)\\
\gvph^{\dagger}_{1...n}(p_1,...,p_n)&=&
(-1)^{\sum_{k=1}^{n-1}([i_1]+...+[i_k])[i_{k+1}]}\,
\bvarphi~^{i_1...i_n}\, (e^{\dagger}_{i_1}\otimes
e^{\dagger}_{i_2}\otimes...\otimes e^{\dagger}_{i_n})\quad
\end{eqnarray}
is a (hermitian) scalar product.

\null

We introduce the finite particle space $\cf_R^0\subset\cf_R$, spanned by
the sequences $(\varphi,\gvph_1,...,$ $\gvph_{1...n},...)$ with
$\gvph_{1...n}\in\ch_R^n$ and $\gvph_{1...n}=0$ for $n$ large
enough. As (\ref{scalar-product}) is defined for all $n$, it
extends naturally to $\cf_R^0$. In this context, the vacuum state
is $\Omega=(1,0,...,0,...)$, so that it is normalized to $1$.

We are now able to define the (smeared) creation and annihilation
operators $A(\gf)$ and $\Adag(\gf)$ on $\cf_R^0$ through their
action: $A(\gf)\Omega=0$ and for $\gvph_{0...n}\in\ch_R^{n+1}$,
\begin{equation}
\label{def-a}
[A(\gf)\gvph]_{1...n}(p_1,...,p_n)=\sqrt{n+1}\int_{\RR}
dp_0\,\gf^{\dagger}_0(p_0)\,\gvph_{0...n}(p_0,p_1,...,p_n)
\end{equation}
Similarly, for $\gvph_{1...n}\in\ch_R^{n}$:
\begin{eqnarray}
\label{def-adag}
&&[\Adag(\gf)\gvph]_{0...n}(p_0,...,p_n)=\frac{1}{\sqrt{n+1}}
\gvph_{1...n}(p_1,...p_n)f_0(p_0)\\
&&{+\frac{1}{\sqrt{n+1}}\sum_{k=1}^n
R_{k-1,k}(p_{k-1}-p_k)...R_{0k}(p_0-p_k)
\gvph_{0...\widehat{k}...n}(p_0,...,\widehat{p_k},...,p_n)\gf_k(p_k)}\nonumber
\end{eqnarray}
where the hatted symbols are omitted.

It is easily checked that (\ref{def-a}) and (\ref{def-adag}) are
indeed elements of $\ch_R^n$ and $\ch_R^{n+1}$ respectively.
Therefore, we have operators acting on $\cf_R^0$ (linearity in
$\gvph$ obvious) with the additional property that they are
bounded (\ie continuous) on each finite particle sector $\ch_R^n$.
Another essential feature is the
adjointness of these operators with respect to $\langle~,~\rangle$
\begin{equation}
\label{adjoint}
\forall\gvph\in\ch_R^n,~\forall\gps\in\ch_R^{n+1},~\forall
\gf\in\ch_R^1,~~~ \langle
\gvph,A(\gf)\gps\rangle=\langle\Adag(\gf)\gvph,\gps\rangle
\end{equation}
At this stage, the Fock representations $\gA(p)$, $\gAdag(p)$ of
the generators of the ZF algebra appear as operator-valued
distributions through the definition
\begin{equation}
A(\gf)=\int_{\RR}dp\, \gf^\dagger(p)\gA(p),~~~~
\Adag(\gf)=\int_{\RR}dp \,\gAdag(p)\gf(p)
\end{equation}
 It is
readily shown from these definitions that $\gA(p)$ and $\gAdag(p)$
satisfy the exchange relations
(\ref{ZF-algebra1}-\ref{ZF-algebra3}) thus providing the desired
representation.

We now have all the ingredients to deduce results for the whole
{F}ock space $\cf_R$ while working on smaller and more intuitive
spaces dense in $\cf_R$, using the continuity of the operators.
In our case, one has to define such a "state space" $\cd\subset\cf_R$
{in the sense of distributions} as follows:
$\cd^0 =\CC$ and
$$
\cd^n =
\left\{\int_{\RR^n}d^np~\gAdag_1(p_1)...\gAdag_n(p_n)\Omega
\gf(p_1,...,p_n);~\gf\in K^nL^2(\RR^n)\right\},\qquad n\ge1
\label{d-space}
$$
Then, $\cd$ is spanned by the sequences
$\gch=(\chi,\gch_1,...,\gch_{1...n},...)$, where
$\gch_{1...n}\in\cd^n$ and $\gch_{1...n}=0$ for $n$ large enough.
We also define
\begin{equation}
\cd_0^0=\CC,~~\cd_0^n=\left\{\tAdag_1(\gf_1,t)...\tAdag_n(\gf_n,t)\Omega,
~\gf_1\succ...\succ \gf_n\right\}\subset\ch_R^n,~n\ge 1
\end{equation}
where
\begin{equation}
\begin{array}{l}
\tAdag(\gf,t)=\int_{\RR}dx~\tgAdag(x,t)\gf(x)\\[1.2ex]
\tgAdag(x,t)=\int_{\RR}dp~\gAdag(p)e^{iqx-iq^2t}
\end{array}\quad x,t\in\RR
\end{equation}
and the space $\cd_0$ is the linear span of sequences
$\gch=(\chi,\gch_1,...,\gch_{1...n},...)$, where
$\gch_{1...n}\in\cd_0^n$ and $\gch_{1...n}=0$ for $n$ large
enough. We also introduce the following partial ordering relation
$$
\gf\succ \bfg~\Leftrightarrow~\forall i,j=1,\ldots,{K},~\forall x\in
supp(f_{i}),~\forall y\in supp(g_{j}),~x>y
$$
which is just the
extension of the ordering of the momenta $k_i$ in the definition
of a state space basis $|k_1,...,k_n\rangle$.
Then, one shows that $\cd$ and $\cd_0$ are dense in $\cf_R$.

Summarizing, we have constructed a graded ZF algebra and its
{F}ock representation $\cf_R$ and, inspired by earlier works
\cite{solQ,skly79,TWC,Dav,DavGut}, we shall see that this allows to
construct the quantum version of NLSS and its solution.

\subsubsection{Quantization of the fields}

{F}ollowing \cite{Dav} and \cite{DavGut}, we simply write the
quantum version of $\phi^{(n)}_j(x,t)$ as
\begin{eqnarray}
\phi^{(n)}_j(x,t)&=&\int_{\RR^{2n+1}}d^n\pp d^{n+1}\qq
\,\sum_{k_{1},\ldots,k_{n}=1}^{K}
\adag_{k_{1}}(p_1)\cdots\adag_{k_{n}}(p_n)a_{k_{n}}(q_n)\cdots
a_{k_{1}}(q_1)a_j(q_0)\nonu &
&\times\frac{e^{i\Omega_{n}(x,t;\pp,\qq)}}{Q_{n}(\pp,\qq,\varepsilon)}
\label{quantum-field-comp}
\end{eqnarray}
using the same notations as in (\ref{solution-rosales2}) and an
$i\eps$ contour prescription. The global field reads
\begin{equation}
\label{quantum-field} \Phi(x,t)=\sum_{n=0}^\infty (-g)^n
\Phi^{(n)}(x,t)~~\text{with}~~\Phi^{(n)}(x,t)=\phi^{(n)}_{j}(x,t)e_j
\end{equation}
{F}rom (\ref{adjoint}), we deduce
\begin{eqnarray}
\Phidag(x,t)&=&\sum_{n=0}^\infty (-g)^n \Phi^{\dagger(n)}(x,t)\\
\text{with}~~\Phi^{\dagger(n)}(x,t)&=&\int_{\RR^{2n+1}}d^n\pp
d^{n+1}\qq
~\gAdag(q_0)\gAdag_1(q_1)\cdots\gAdag_n(q_n)\gA_n(p_n)\cdots
\gA_1(p_1)\nonu &
&\times\frac{e^{-i\Omega_{n}(x,t;\pp,\qq)}}{Q_{n}(\pp,\qq,-\varepsilon)}
\end{eqnarray}
Just like we dealt with $A(\gf)$ and $\Adag(\gf)$, we are
 naturally led to introduce
\begin{equation}
\Phi(\gf,t)=\int_{\RR}\gf^\dagger(x)\Phi(x,t),~~\Phidag(\gf,t)=
\int_{\RR}\Phidag(x,t)\gf(x)
\end{equation}
And just like we did in \cite{CR}, one shows that $\Phi(\gf,t)$
and $\Phidag(\gf,t)$ are indeed well-defined operators on a common
invariant domain which turns out to be $\cd_0$. These fields also
satisfy the following fundamental requirement

\begin{theo}
The quantum fields $\Phi(\gf,t),~\Phidag(\bfg,t)$ satisfy the
equal time canonical commutation relations as operators on
$\cf_R^0$
\begin{eqnarray}
\label{global-CCR1}
[\Phi(\gf,t),\Phi(\bfg,t)] &=& [\Phidag(\gf,t),\Phidag(\bfg,t)]=0\\
\label{global-CCR2} [\Phi(\gf,t),\Phidag(\bfg,t)] &=& \langle
\gf,\bfg\rangle
\end{eqnarray}
\end{theo}
\prf the proof is the same as in the ordinary NLS equation, see
\cite{Dav} or \cite{GLM} for details.
\finprf

One then deduces the equal time CCR in
components for the operator-valued distributions
$\phi_j(x,t),~\phidag_k(y,t)$:
\begin{eqnarray}
\label{comp-CCR1}
[\hspace{-2pt}[\phi_j(x,t),\phi_k(y,t)]\hspace{-2pt}] &=&
[\hspace{-2pt}[\phidag_j(x,t),\phidag_k(y,t)]\hspace{-2pt}]=0
\\
\label{comp-CCR2}
[\hspace{-2pt}[\phi_j(x,t),\phidag_k(y,t)]\hspace{-2pt}] &=&
\delta_{jk}\delta(x-y)
\end{eqnarray}
Let us remind that for $j,k=M+1,\ldots K$, the above CCR correspond to
anticommutator, consistent with the fermionic nature of these fields.

\subsubsection{Time evolution\label{timeevo}}

We first wish to emphasize that the form of the Hamiltonian
(\ref{defi-hamiltonian}) cannot be reproduced here owing to the
nature of the fields (products of distributions are not defined).
{F}ortunately, the power of the ZF algebra and the quantum inverse
method (leading to (\ref{quantum-field-comp}-\ref{quantum-field}))
rescues us by delivering a simple, free-like Hamiltonian in terms
of oscillators. Indeed, one easily checks that the Hamiltonian
defined by
\begin{equation}
H=\int_{\RR}dp~p^2\gAdag(p)\gA(p)
\end{equation}
is self-adjoint, \ie $H^\dagger=H$. Moreover,
\begin{equation}
\label{eigen} \forall\gvph\in\cd,~~[H\gvph]_{1...n}(p_1,...,p_n)=
(p_1^2+...+p_n^2)\gvph_{1...n}(p_1,...,p_n)
\end{equation}
which shows that $\cd$ is also an invariant domain for $H$ and
that this operator has the correct eigenvalues. {F}inally, $H$
generates the time evolution of the field:
\begin{equation}
\label{time} \Phi(f,t)=e^{iHt}\Phi(f,0)e^{-iHt}
\end{equation}
Therefore, $H$, so defined, is the Hamiltonian of our quantum
system.

Note that (\ref{eigen}) and (\ref{time}) have to be understood as
operator equalities and must be evaluated on $\cd$.

The free-like expression for $H$ in terms of creation and
annihilation oscillators may be surprising at first glance but it
is actually a mere consequence of the rather complicated exchange
relations (\ref{ZF-algebra1}-\ref{ZF-algebra3}). One can say that
the effect of the non-linear term has been encoded directly in the
oscillators instead of the Hamiltonian (or equivalently the
Lagrangian) of the field theory, yielding a (possibly misleading)
simple expression for $H$. One may finally wonder about the
coupling constant which seems to disappear. Once again, it is
actually present through the $R$-matrix in the exchange relations.

 Besides, the quantum
nonlinear super-Schr\"odinger equation holds in the following
form:
\begin{equation}
\forall
\gvph,\gps\in\cd,~~(i\partial_t+\partial_x^2)\langle\gvph,\Phi(x,t)\gps\rangle
=2g\langle\gvph,:\Phi\Phidag\Phi:(x,t)\gps\rangle
\end{equation}

\subsubsection{Correlation functions}

Again following the case of NLS, one shows that for
$\gvph,\gps\in\cd$, one has
\begin{equation}
\label{psidag1}
\gf\succ \bfg,\quad \langle\gvph,\Phidag(\bfg,t)\tAdag(\gf,t)\gps\rangle=
\langle\gvph,\tAdag(\gf,t)\Phidag(\bfg,t)\gps\rangle
\end{equation}
 for $\bfg\succ \gf_i,~i=1,...,n$
\begin{equation}
\label{psidag2}
\langle\gvph,\Phidag(\bfg,t)\tAdag(\gf_1,t)...\tAdag(\gf_n,t)\Omega\rangle=
\langle\gvph,\tAdag(\bfg,t)\tAdag(\gf_1,t)...\tAdag(\gf_n,t)\Omega\rangle
\end{equation}
 for any $\gf_1\succ \gf_2\succ ...\succ \gf_n$
\begin{equation}
\hspace{-2.1ex}
\langle\gvph,\Phi(\bfg,t)\tAdag(\gf_1,t)...\tAdag(\gf_n,t)\Omega\rangle=
\sum_{j=1}^n\langle
\bfg,\gf_j\rangle\langle\gvph,\tAdag(f_1,t)...\widehat{\tAdag}(\gf_j,t)...
\tAdag(\gf_n,t)\Omega\rangle\label{psi3}
\end{equation}
This proves that the correlation functions of the NLSS model are
completely determined, e.g.:
$$
\langle\Omega,\Phi(\bfg_{1},t)\cdots\Phi(\bfg_{m},t)
\Phidag(\gf_1,t)\cdots\Phidag(\gf_n,t)\Omega\rangle=
\delta_{m,n}\sum_{\sigma\in S_{n}}\prod_{i=1}^n\langle
\bfg_{\sigma(i)},\gf_{i}\rangle
$$
\begin{eqnarray*}
\lefteqn{\langle\gvph_{1\ldots p},\Phi(\bfg_{1},t)\cdots\Phi(\bfg_{n},t)
\Phidag(\gf_1,t)\cdots\Phidag(\gf_m,t)\Omega\rangle=}
\qquad\qquad\\
&&\qquad\qquad
\delta_{m,n+p}\sum_{\sigma\in S_{n+p}}\left(\prod_{i=1}^n\langle
\bfg_{\sigma(i)},\gf_{i}\rangle\right)
\langle \gvph_{1\ldots p},\bfg_{\sigma(n+1)}\cdots\bfg_{\sigma(n+p)}\rangle
\end{eqnarray*}
Similar expressions can be obtained when dealing with the fields $\Phi(x,t)$ and $\Phidag(x,t)$.
\section{Lax pair and super-Yangian symmetry for NLSS\label{sec3}}

Let us stress once again that we aim at generalizing known results
of integrability and symmetry for the non-linear Schr\"odinger
equation to the case of an arbitrary number of bosons \textbf{and}
fermions. This physical motivation can be carried out by using
appropriately the graded formalism presented in the appendix.
Furthermore, we also want to transport our results to the quantum
case, which leads us to adopt the convenient Hamiltonian form of
our model.

\subsection{Classical Lax pairs \label{Lax}}
We define the Lax even super-matrix in $gl(M+1|N)$

\begin{eqnarray}
L(\lda;x) &=& \frac{i\lda}{2}\Sigma+\Omega(x)\mb{with}
\Sigma=\II_{K+1,K+1}-2E_{K+1,K+1} \label{Laxmatrix1}\\
\mb{and}&&
\Omega(x)=i\sqrt{g}\,\sum_{j=1}^{K}\Big(\phi_{j}(x)E_{j,K+1}
-\phidag_{j}(x)E_{K+1,j}\Big) \label{Laxmatrix2}
\end{eqnarray}

Let us stress that, as above, the elementary matrices $E_{jk}$
(with 1 at position $j,k$) are $\ZZ_{2}$-graded, with
$[E_{jk}]=[j]+[k]$, $[j]=[K+1]=0$ for $1\le j\le M$ and $[j]=1$
for $M<j\le K$. With this convention, the $gl(M+1|N)$ superalgebra
has the unusual matrix form
$$\left(\begin{array}{c|c|c}
\rule{0ex}{3.5ex}\ M\times M\  &  &\! M\!\times \!1\! \\[1ex] \hline
 \rule{0ex}{3.5ex}& \ N\times N\  & \\[1ex] \hline 1\times M &  &
\! \!1\times 1\! \end{array}\right)$$
where the size of the sub-matrices corresponding to
bosonic generators have been explicitly written.

Using the PB of the $\phi$'s, it is easy to compute that
\begin{eqnarray}
\{L_{1}(\lda;x), L_{2}(\mu;y)\} &=& i\delta(x-y)\,
\left[r(\lda-\mu), L_{1}(\lda;x)+ L_{2}(\mu;y)\right]
\\
\mb{with} && r(\lda-\mu) = \frac{g}{\lda-\mu}\ \Pi_{12}
\label{defr}
\end{eqnarray}
where we have introduced the $(K+1)\times (K+1)$ super-permutation
$$
\Pi_{12}=\sum_{i,j=1}^{K+1}\, (-1)^{[j]}\, E_{ij}\otimes E_{ji}
$$

\begin{defi}
    We define the transition matrix by
\begin{equation}
\prt_{x}T(\lda;x,y) = L(\lda;x) T(\lda;x,y),\ x>y
\end{equation}
with the "initial condition" $T(\lda;x,x)=\II$.
\end{defi}
$T(\lda;x,y)$ obeys the iterative equation
\begin{equation}
T(\lda;x,y) = E(\lda;x-y)+E(\lda;x)\int_{y}^x dz\,
\Omega(z)E(\lda;z)T(\lda;z,y)\label{ite-1}
\end{equation}
where we have introduced
\begin{equation}
    E(\lda;x)=\exp\left(\frac{ix\lda}{2}\Sigma\right)=
    e^{ix\lda/2}\,\II_{K+1}+\left(e^{-ix\lda/2}-
    e^{ix\lda/2}\right)E_{K+1,K+1}
\end{equation}
\begin{prop}
     \begin{equation}
\{T_{1}(\lda;x,y), T_{2}(\mu;x,y)\} = \left[r(\lda-\mu),
T(\lda;x,y)\otimes T(\mu;x,y)\right] \label{comTxTy}
     \end{equation}
\end{prop}
\prf The eq. (\ref{ite-1}) implies that
\begin{eqnarray}
T(\lda;x,y) &=& \sum_{n=0}^\infty T^{(n)}(\lda;x,y)\\
T^{(n)}(\lda;x,y) &=& \int_{\RR^n} d^nz\,
\theta(x>z_{1}>z_{2}>\ldots>z_{n}>y)\,
E(\lda;x-z_{1})\Omega(z_{1})\nonu &&\times
E(\lda;z_{1}-z_{2})\Omega(z_{2})\cdots
\Omega(z_{n})E(\lda;z_{n}-y)\qquad \label{Tn-1}
\end{eqnarray}
It is then simple to show that
\begin{eqnarray}
\left\{\Phi_{1}(w),T_{2}(\lda;x,y)\right\} &=& \sqrt{g}\,
\theta(x>w>y) \,T_{2}(\lda;x,w)\,\sigma_{12}^-\,T_{2}(\lda;w,y)
\label{PBphiT}\\
\left\{\Phi_{2}(w),T_{1}(\lda;x,y)\right\} &=& \sqrt{g}\,
\theta(x>w>y) \,T_{1}(\lda;x,w)\,\sigma_{21}^-\,T_{1}(\lda;w,y)
\label{PBTphi}\\
\left\{\Phidag_{1}(w),T_{2}(\lda;x,y)\right\} &=& \sqrt{g}\,
\theta(x>w>y)\, T_{2}(\lda;x,w)\,\sigma_{12}^+\,T_{2}(\lda;w,y)
\label{PBbphT}\\
\left\{\Phidag_{2}(w),T_{1}(\lda;x,y)\right\} &=& \sqrt{g}\,
\theta(x>w>y)\, T_{1}(\lda;x,w)\,\sigma_{21}^+\,T_{1}(\lda;w,y)
\label{PBTbph}
\end{eqnarray}
where we have defined
\begin{equation}
\! \sigma_{12}^- = \sum_{j=1}^{K} e_{j}\otimes E_{K+1,j}\  ;\
\sigma_{12}^+ =\sum_{j=1}^{K} (-1)^{[j]}\, e^\dagger_{j}\otimes
E_{j,K+1}\
\end{equation}
{F}rom the form (\ref{PB-sg})
%and the identity (\ref{deltaphi}),
one also computes
\begin{eqnarray}
\left\{\Phi_{1}(w),T_{2}(\lda;x,y)\right\} &=&
i\,(e_{j}\otimes\II)\,\frac{\delta\,T_{2}(\lda;x,y)}{\delta\,\phidag_{j}(w)}\\
\left\{\Phi_{2}(w),T_{1}(\lda;x,y)\right\} &=&
i\,(\II\otimes e_{j})\,\frac{\delta\,T_{1}(\lda;x,y)}{\delta\,\phidag_{j}(w)}\\
\left\{\Phidag_{1}(w),T_{2}(\lda;x,y)\right\} &=&
-i\,(-1)^{[j]}(e_{j}^\dagger\otimes\II)\,\frac{\delta\,T_{2}(\lda;x,y)}{\delta\,\phi_{j}(w)}\\
\left\{\Phidag_{2}(w),T_{1}(\lda;x,y)\right\} &=&
-i\,(-1)^{[j]}(\II\otimes
e_{j}^\dagger)\,\frac{\delta\,T_{1}(\lda;x,y)}{\delta\,\phi_{j}(w)}
\end{eqnarray}
This shows that the PB can be rewritten as
\begin{eqnarray}
\hspace{-2.4ex}
    \left\{T_{1}(\lda,x,y),T_{3}(\mu;x,y)\right\} &=&
i\int_{\RR}dw\,\Big(
\{\Phidag_{2}(w),T_{1}(\lda;x,y)\}\,\{\Phi_{2}(w),T_{3}(\mu;x,y)\}\nonu
&&\quad -
\{\Phidag_{2}(w),T_{3}(\mu;x,y)\}\,\{\Phi_{2}(w),T_{1}(\lda;x,y)\}\Big)
\end{eqnarray}
Inserting (\ref{PBphiT}) and (\ref{PBbphT}) in this expression,
 one gets
 \begin{eqnarray}
\left\{T_{1}(\lda;x,y),T_{2}(\mu;x,y)\right\} &=& ig
\int_{y}^xdw\, T_{1}(\lda;x,w)T_{2}(\mu;x,w)\, (\pi_{12}-\pi_{21})\,
T_{1}(\lda;w,y)T_{2}(\mu;w,y)
\nonu
\mb{where}\pi_{12} &=& \sum_{j=1}^{K}\, E_{j,K+1}\otimes
E_{K+1,j}\label{eq:pi12}
\end{eqnarray}
Finally, a direct calculation shows that
\begin{equation}\begin{array}{l}
\displaystyle \frac{\prt}{\prt
w}\Big(T_{1}(\lda;x,w)T_{2}(\mu;x,w)\,\Pi_{12}\,T_{1}(\mu;w,y)T_{2}(\lda;w,y)\Big) \\
\displaystyle\qquad=
i\,\frac{\lda-\mu}{2}\,T_{1}(\lda;x,w)T_{2}(\mu;x,w)\,
(\pi_{12}-\pi_{21})\,
T_{1}(\lda;w,y)T_{2}(\mu;w,y)
\end{array}
\end{equation}
 so that we get (\ref{comTxTy}).
\finprf
\begin{prop}
    The following limits are well-defined:
\begin{eqnarray}
T^-(\lda;x)&=&\lim_{y\to-\infty}T(\lda;x,y)E(\lda;y)\label{defT-}\\
T^+(\lda;y)&=&\lim_{x\to\infty}E(\lda;-x)T(\lda;x,y)\label{defT+}\\
T(\lda)&=&T^+(\lda;z)T^-(\lda;z) \ =\ \lim_{x\to\infty\atop
y\to-\infty}E(\lda;-x)T(\lda;x,y)E(\lda;y) \label{defTmono}
\end{eqnarray}
$T(\lda)$ is called the monodromy
matrix.
\end{prop}
\prf Using the equality
$E(\lda;x)\Omega(z) = \Omega(z)E(\lda;-x)$,
valid for any $x$, $z$,
$T^{(n)}(\lda;x,y)$ can be conveniently rewritten as
\begin{eqnarray}
T^{(n)}(\lda;x,y) &=& E(\lda;x)\int_{\RR^n} d^nz\,
\theta(x>z_{1}>\ldots>z_{n}>y)\nonu &&\times\,E(\lda;2\sum_{j=1}^n
(-1)^jz_{j}) \left(\prod_{k=1}^n\Omega(z_{k})\right)E(\lda;-y)
\label{Tn-2}
\end{eqnarray}
which shows that the limits are well-defined. \finprf

\begin{prop}
\begin{equation}
\label{classicalRTT}
 \{T_{1}(\lda), T_{2}(\mu)\} = r_{+}(\lda-\mu) T(\lda)\otimes
T(\mu) - T(\lda)\otimes T(\mu) r_{-}(\lda-\mu)
\end{equation}
with
\begin{eqnarray}
    r_{+}(\lda-\mu) &=& \frac{g}{\lda-\mu}\,\left( P_{12}+
    E_{K+1,K+1}\otimes E_{K+1,K+1}\right)\\
    &&+i\pi g\delta(\lda-\mu)\,(\pi_{12}-\pi_{21})
    \nonu
    r_{-}(\lda-\mu) &=& \frac{g}{\lda-\mu}\left( P_{12}+
    E_{K+1,K+1}\otimes E_{K+1,K+1}\right)\\
    &&-i\pi g\delta(\lda-\mu)\,(\pi_{12}-\pi_{21})
\nonumber
\end{eqnarray}
where $P_{12}$ is the super-permutation in the space of $K\times
K$ matrices.
\end{prop}
\prf Direct calculation, plugging (\ref{defTmono}) into
(\ref{comTxTy}) and using the Cauchy principal value
$\lim_{\lda\to\infty}p.v.\left(\frac{e^{\pm i\lda
x}}{x}\right)=\pm i\pi\delta(x)$. \finprf

Introducing $t(\lda)$, the $K\times K$ sub-matrix of $T(\lda)$ with
the last row and column removed, and $D(\lda)=T_{K+1,K+1}(\lda)$,
one finally computes:
\begin{eqnarray}
    \{t_{1}(\lda),t_{2}(\mu)\} &=& \frac{g}{\lda-\mu}\,[P_{12}\,,\,t(\lda)\otimes
    t(\mu)] \label{SY}\\
    \{D(\lda),t(\mu)\} &=& 0\label{PB-Dt}\\
    \{D(\lda),D(\mu)\} &=& 0\label{DD}
\end{eqnarray}
(\ref{SY}) shows that $t(\lda)$ defines a classical version of the
super-Yangian $Y(gl(M|N))$.
(\ref{DD}) shows that $D(\lda)$ can be taken as a generating
function for a hierarchy, and (\ref{PB-Dt}) proves that the
super-Yangian is a symmetry of this hierarchy. It remains to
identify this hierarchy.
\begin{lem}
    Only $T^{(2n)}(\lda)$, $n\in\ZZ_{+}$, contribute to the super-Yangian
    generators $t(\lda)$ and to the Hamiltonian generating function
    $D(\lda)$. \\
    Expanding $t(\lda)$ and $D(\lda)$ as series in $\lda^{-1}$,
    one has $T^{(2n)}(\lda)=o(\lda^{-n})$.
\end{lem}
\prf It is clear that $T^{(n)}(\lda)$ contains the product of
exactly $n$ matrices $\Omega$, the other matrices entering in its
definition being diagonal. Due to the form of $\Omega$, only
products of an even number of such matrices will contribute to
$t(\lda)$ and $D(\lda)$.

To show the $\lambda$ dependence, we consider the integration on
$z_{2j}$ and $z_{2j+1}$, and perform an integration by part,
assuming that the fields $\Phi$ and $\Phidag$ are vanishing at
infinity:
\begin{eqnarray*}
&&  \int_{-\infty}^{z_{2j-1}} dz_{2j}\, \int_{-\infty}^{z_{2j}}
dz_{2j+1}\, E(\lda;2z_{2j}-2z_{2j+1})\Omega(z_{2j})
\Omega(z_{2j+1}) I_{j,n}(z_{2j+1},\ldots,z_{2n})\ =\\
&& \frac{i}{\lda}\Sigma \int_{-\infty}^{z_{2j-1}} dz_{2j}\,\left[
 \Omega(z_{2i})^2 I_{j,n}(z_{2j+1},\ldots,z_{2n})-\int_{-\infty}^{z_{2j}} dz_{2j+1}\,
 E(\lda;2z_{2j}-2z_{2j+1})
 \times\right.\\
&&\left. \times\Omega(z_{2j})\prt_{2j+1}\Big( \Omega(z_{2j+1})
I_{j,n}(z_{2j+1},z_{2j+2},\ldots,z_{2n})\Big) \right]
\end{eqnarray*}
Above, $\prt_{k},\ \forall k$ stands for $\frac{\prt}{\prt
z_{k}}$, and $I_{j,n}(z_{2j+1},z_{2j+2},\ldots,z_{2n})$ denotes
the other integrals (depending on $z_{k}$, $k\geq 2j$) which
enters into the definition of $T^{(n)}(\lda)$.

It is clear that one can do this integration for all  $z_{2j}$,
$j=1,\ldots,n$ and any number of times, so that the lowest power
of $\lda^{-1}$ is $n$. \finprf
\begin{prop}\label{firsthami}
The first Hamiltonians generated by $D(\lda)$ read
\begin{eqnarray}
D^{(1)} &=&igN\mb{with} N = \int_{-\infty}^\infty dx\, \Phidag(x)\,\Phi(x) \\
D^{(2)} &=& -\half g^2N^2+gP \mb{with} P= \int_{-\infty}^\infty
dx\,
\Phidag(x)\,\prt\Phi(x) \\
D^{(3)} &=& -\frac{ig^3}{6}N^3+ig^2\,NP +igH\\
H &=& \int_{-\infty}^\infty dx\, \prt\Phidag(x)\,\prt\Phi(x)
+g\int_{-\infty}^\infty dx\, \left(\Phidag(x)\,\Phi(x)\right)^2
\end{eqnarray}
This shows that $D(\lda)$ generates the Hamiltonians of the NLSS
hierarchy, so that (\ref{PB-Dt}) proves that $Y(gl(M|N))$ is a
symmetry of this hierarchy.
\end{prop}
\prf We use the technics given in the above proof, focusing on the
$(K+1,K+1)$ matrix element. The bound in the integrals are
simplified using the property
\begin{equation}\begin{array}{l}
\Big(\Omega(x_{1})\prt^k\Omega(x_{2})\Omega(x_{3})\prt^l\Omega(x_{4})\Big)_{K+1,K+1}
=\Big(\Omega(x_{1})\prt^k\Omega(x_{2})\Big)_{K+1,K+1}\times\\
\hspace{17em}\times\Big(\Omega(x_{3})\prt^l\Omega(x_{4})\Big)_{K+1,K+1}
\end{array}
\end{equation}
\finprf

\subsection{Quantum Lax pair}

Following Sklyanin \cite{Skly}, we define
\begin{defi}
The quantum transition matrix $\ct(\lda;x,y)$ is the Wick
(normal)-ordered classical transition matrix $T(\lda;x,y)$
regarded as a functional of the quantum canonical fields
$\Phi(x),\Phidag(x)$
\begin{equation}
\ct(\lda;x,y)=~:T(\lda;x,y):
\end{equation}
Here and below the normal ordering is defined as
$$
:\phi_{j}(x)\phidag_{k}(y):=(-1)^{[j][k]}\phidag_{k}(y)\phi_{j}(x),\quad
\forall x,y
$$
and extended to monomials in $\phi$, $\phidag$ in the usual way, \ie
with all the $\phi$'s on the right of the $\phidag$'s, keeping the
original order between the $\phi$'s and between the $\phidag$'s.
\end{defi}

For convenience,
we also define a symbol $\ddagger~\ddagger$ which acts on operators and is not to be
confused with the symbol $:~:$. It simply guarantees the ordering
of $\Phi, \Phidag$ in an expression containing $L(\lda;x)$ and other
(normal-ordered) functionals
of the quantum fields without changing the internal ordering of
the functionals. For example, if $A=:a:$ and $B=:b:$ then
\begin{eqnarray*}
    \ddagger AL(\lda;x) B\ddagger &=&\frac{i\lda}{2}A\Sigma B\\
    &&+i\sqrt{g}\sum_{j=1}^{K}\left(
(-1)^{[j][A]}\phi_j(x)AE_{K+1,j}B-(-1)^{[j][B]}AE_{j,K+1}B\phidag_j(x)\right)
\end{eqnarray*}

The previous definition
gives rise to many questions dealing with operator theory and
functional analysis which were answered for the bosonic case in
the very detailed review \cite{Gut} by Gutkin. But for the sake of
brevity, we mimic the
compact, albeit more formal, approach of Sklyanin since it
contains all the fundamental and physical ideas, bearing in mind
that everything is well-defined.

In this sense, the quantum transition matrix
is the fundamental solution of the quantum auxiliary problem
\begin{equation}
\label{defLT}
\prt_{x}\ct(\lda;x,y)=\ddagger L(\lda;x)\ct(\lda;x,y)\ddagger
 ~~\text{with}~~\ct(\lda;x,x)=\1
\end{equation}
and satisfies
\begin{eqnarray}
\prt_{y}\ct(\lda;x,y)&=&-\ddagger \ct(\lda;x,y)L(\lda;y)\ddagger\nonu
 \ct(\lda;x,y)\ct(\lda;y,z) &=&\ct(\lda;x,z)~~for~~x<y<z~~or~~x>y>z\nonumber
\end{eqnarray}
where $L(\lda;x)$ is the Lax even super-matrix defined in
(\ref{Laxmatrix1})-(\ref{Laxmatrix2}).

This system of first order differential equations together with
the given initial condition  is equivalent to the following
Volterra integral representations:
\begin{eqnarray}
\ct(\lda;x,y)&=&\1 +\int_y^x d\omega\ddagger
L(\lda;\omega)\ct(\lda;\omega,y)\ddagger \label{volt1}\\
\ct(\lda;x,y)&=&\1 +\int_y^x d\omega\ddagger
\ct(\lda;x,\omega)L(\lda;\omega)\ddagger \label{volt2}
\end{eqnarray}

In order to reach our final goal there are several steps which all
rely on one simple idea extensively used in the inverse problem
literature, that is two quantities are equal if and only if they
satisfy the same first order differential equation with the same
initial condition. This is what is called "the differential
equation approach" by Gutkin in \cite{Gut}. He criticized  this
approach but showed that it gives the correct answer using the
"discrete approximation approach" which amounts to the same line
of argument but deals with finite differences on subintervals
 of $[x,y]$ instead of a true derivative.

The first step is to obtain the commutation relations of matrix
elements of the transition matrix and we need two preliminary
lemmas.

\begin{lem} \label{deriveTT}
$\ct_1(\lda;x,y)\ct_2(\mu;x,y)$ satisfies the following differential
system:
\begin{eqnarray}
\prt_x \{\ct_1(\lda;x,y)\ct_2(\mu;x,y)\}&=&\ddagger
\cl_{12}(\lda,\mu;x)\ct_1(\lda;x,y)\ct_2(\mu;x,y)\ddagger\label{deriveTT1}\\
\ct_1(\lda;x,x)\ct_2(\mu;x,x)&=&\ct_2(\mu;x,x)\ct_1(\lda;x,x)=\1\otimes\1
\end{eqnarray}
where
\begin{equation}
\cl_{12}(\lda,\mu;x)=L_1(\lda;x)+L_2(\mu;x)+g\,\pi_{12}
\end{equation}
\end{lem}

\prf The idea is once again to use the equivalence between the
differential problem and the Volterra integral representation of
the solution. Indeed, taking care of the ordering of the fields
when using (\ref{volt1},\ref{volt2}), one gets
\begin{equation*}
\ct_1(\lda;x,y)\ct_2(\mu;x,y)=\1\otimes\1 +\int_y^x
d\omega\ddagger\cl_{12}(\lda,\mu;\omega)\ct_1(\lda;\omega,y)\ct_2(\mu;\omega,y)\ddagger
\end{equation*}
\finprf

\begin{lem}\label{RL}
The operator $\cl_{12}(\lda,\mu;x)$  satisfy the
following relation:
\begin{equation}
\calr_{12}(\lda-\mu)\cl_{12}(\lda,\mu;x)=\cl_{21}(\mu,\lda;x)\calr_{12}(\lda-\mu)
\end{equation}
where $\calr_{12}(\lda-\mu)=\1 -ir(\lda -\mu)$, and $r(\lda -\mu)$
is given by (\ref{defr}).
\end{lem}

\prf  direct calculation using
$$
\left[\Pi_{12},L_1(\lda;x)+L_2(\mu;x)\right]=
i(\lda-\mu)\,(\pi_{12}-\pi_{21})
$$
where $\pi_{12}$ has been defined in (\ref{eq:pi12}).
 \finprf

We can now formulate the basic result of this paragraph.
\begin{theo}\label{finitevolRTTtheo}
The quantum transition matrix $\ct(\lda;x,y)$ satisfies the
following finite volume commutation relations:
\begin{equation}
\label{finitevolRTT}
\calr_{12}(\lda-\mu)\ct_1(\lda;x,y)\ct_2(\mu;x,y)=\ct_2(\mu;x,y)\ct_1(\lda;x,y)\calr_{12}(\lda-\mu)
\end{equation}
\end{theo}

\prf Using the fact that $\calr_{12}(\lda)$ is a numerical,
invertible (for $\lda$ real and nonzero) matrix, lemmas
\ref{deriveTT} and \ref{RL} imply that the quantities
$\ct_2(\mu;x,y)\ct_1(\lda;x,y)$ and
$\calr_{12}(\lda-\mu)\ct_1(\lda;x,y)\ct_2(\mu;x,y)\calr_{12}^{-1}(\lda-\mu)$
satisfy the same first order differential equation with the same
initial condition. \finprf

Let us remark that if we restore the Planck constant in the
canonical commutation relations, then $\calr_{12}(\lda-\mu)=\1
-i\hbar r(\lda -\mu)$ and we recover the relation (\ref{comTxTy})
for the classical transition matrix, given that as
$\hbar\rightarrow 0$, $\ct(\lda;x,y)\rightarrow T(\lda;x,y)$ and
$[~,~]\rightarrow i\hbar\{~,~\}$ and keeping the terms of order
$\hbar$.

 We are now
in position to define the quantum monodromy matrix as an
appropriate limit of the quantum transition matrix to obtain the
infinite volume commutation relations corresponding to
(\ref{finitevolRTT}). The
crucial difference with respect to the classical case comes from
the nontrivial commutation relations of the quantum fields, which
produces the term proportional to $g$ in $\cl_{12}(\lda,\mu;x)$.

Therefore, one cannot define the limit as in (\ref{defTmono}) and
insert it directly in the finite volume commutation relations.
Instead, we are led to compare the asymptotic behaviour of
$\ct_1(\lda;x,y)\ct_2(\mu;x,y)$, for which we have information
with that of $\ct_1(\lda;x,y),\ct_2(\mu;x,y)$ separately, whose
commutation relations in the infinite interval limit we are
looking for.

\begin{defi}
The quantum equivalents of (\ref{defT-}-\ref{defTmono}) are
defined by
\begin{equation}
\ct^-(\lda;x)=\,:T^-(\lda;x):,~~\ct^+(\lda;y)=\,:T^+(\lda;y):,~~\ct(\lda)=\,:T(\lda):
\end{equation}
and $\ct(\lda)=\ct^+(\lda;z) \ct^-(\lda;z)$ is the quantum
monodromy matrix.
\end{defi}
$E(\lda;x)$ being a numerical matrix, one immediately deduces
\begin{eqnarray}
\prt_x\ct^-(\lda;x)&=&\ddagger L(\lda;x)\ct^-(\lda;x)\ddagger\label{deriveT-}\\
\prt_x\ct^+(\lda;x)&=&-\ddagger\ct^+(\lda;x)L(\lda;x)\ddagger
\end{eqnarray}
As a first step, we look for information on
 $\ct^-_1(\lda;x)\ct^-_2(\mu;x)$ from what
we know of $\ct_1(\lda;x,y)\ct_2(\mu;x,y)$. This is gathered in

\begin{lem}\label{lemlambda}

\begin{equation}
\lim_{y\to
-\infty}\ct_1(\lda;x,y)\,\ct_2(\mu;x,y)\,\xi_{12}(\lda,\mu;y)=
\ct^-_1(\lda;x)\,\ct^-_2(\mu;x)C_{12}(\lda,\mu)
\end{equation}
where, $\pi_{12}$ being defined as in (\ref{eq:pi12}), we have
introduced:
\begin{eqnarray}
\xi_{12}(\lda,\mu;y)&=&\exp
\left[(\frac{i\lda}{2}\Sigma_1+\frac{i\mu}{2}\Sigma_2
+g\,\pi_{12})y\right]\\
C_{12}(\lda,\mu)&=&\1\otimes\1-\frac{ig}{\lda-\mu+i\varepsilon}\,\pi_{12}\label{defC}
\end{eqnarray}
\end{lem}

\prf Let
\begin{eqnarray}
\Lda(\lda,\mu;x)&=&\lim_{y\to
-\infty}\ct_1(\lda;x,y)\,\ct_2(\mu;x,y)\,\xi_{12}(\lda,\mu;y)\\
\Lda^-(\lda,\mu;x)&=&\ct^-_1(\lda;x)\,\ct^-_2(\mu;x)
\end{eqnarray}
Rewriting
$\cl_{12}(\lda,\mu;x)=\cl_0(\lda,\mu)+\Omega_1(x)+\Omega_2(x)$
with
$\cl_0(\lda,\mu)=\frac{i\lda}{2}\Sigma_1+\frac{i\mu}{2}\Sigma_2+g
\,\pi_{12}$, one easily gets from (\ref{deriveTT1}) the integral
representation
\begin{eqnarray}
&&\ct_1(\lda;x,y)\ct_2(\mu;x,y) =\xi_{12}(\lda,\mu;x-y)
\label{repTT}\\
&&{\qquad\qquad+\int_y^xd\omega\ddagger
\ct_1(\lda;x,\omega)\ct_2(\mu;x,\omega)
\left(\Omega_1(\omega)+\Omega_2(\omega)\right)
\ddagger\xi_{12}(\lda,\mu;\omega-y)}\nonumber
\end{eqnarray}
which shows that $\Lda(\lda,\mu;x)$ is well-defined and also
satisfies
$$\prt_x\Lda(\lda,\mu;x)=\ddagger\cl_{12}(\lda,\mu;x)\Lda(\lda,\mu;x)\ddagger$$
Now following the same line of argument as in lemma
\ref{deriveTT}, we get
$$\prt_x\Lda^-(\lda,\mu;x)=\ddagger\cl_{12}(\lda,\mu;x)\Lda^-(\lda,\mu;x)\ddagger$$
Consequently,
\begin{equation}
\label{differ}
\Lda(\lda,\mu;x)=\Lda^-(\lda,\mu;x)C_{12}(\lda,\mu),~~\forall x
\end{equation}
and we can determine $C_{12}(\lda,\mu)$ from the asymptotic behaviour
as $x\to -\infty$. From the physical requirement that
$$\lim_{x\to\pm\infty}|\Phi(x)|=0$$ and eq.(\ref{repTT}),
we see that
\begin{equation*}
\ct_1(\lda;x,y)\ct_2(\mu;x,y)\underset{x\to y}{\underset{y\to
-\infty}{\sim}}\xi_{12}(\lda,\mu;x-y)
\end{equation*}
implying
\begin{equation}
\Lda(\lda,\mu;x)\underset{x\to-\infty}{\sim}\xi_{12}(\lda,\mu;x)
\end{equation}
On the other hand, from (\ref{deriveT-}), $\Lda^-(\lda,\mu;x)$ can
be represented as
\begin{eqnarray*}
\Lda^-(\lda,\mu;x)&=&E_1(\lda;x)E_2(\mu;x)
+\int_{-\infty}^xd\omega\ddagger\ct_1(\lda;x,\omega)\ct_2(\mu;x,\omega)\\
&&\times\left(\Omega_1(\omega)+\Omega_2(\omega)+g\pi_{12}
\right)\ddagger E_1(\lda;\omega)E_2(\mu;\omega)
\end{eqnarray*}
so that
$$
\Lda^-(\lda,\mu;x)\underset{x\to-\infty}{\sim}E_1(\lda;x)E_2(\mu;x)+I(\lda,\mu;x)
$$
where
$$
I(\lda,\mu;x)=g\int_{-\infty}^xd\omega\,\xi_{12}(\lda,\mu;x-\omega)\,
\pi_{12}\,E_1(\lda;\omega)E_2(\mu;\omega)
$$
can be evaluated from the knowledge of
$$
\xi_{12}(\lda,\mu;x)=E_1(\lda;x)E_2(\mu;x)+2g\,
\frac{sin\left(\frac{\lda-\mu}{2}x\right)}
{\lda-\mu}\,\pi_{12}
$$
and an $i\varepsilon$ prescription to get
$$
I(\lda,\mu;x) = \frac{ig}{\lda-\mu+i\varepsilon}\,
e^{-i\frac{\lda-\mu}{2}x}\,\pi_{12}
$$
Now, adopting the regularization
$$
2g\frac{sin\left(\frac{\lda-\mu}{2}x\right)}{\lda-\mu} =
\frac{-ig}{\lda-\mu+i\varepsilon} \left[e^{i\frac{\lda-\mu}{2}x}
-e^{-i\frac{\lda-\mu}{2}x}\right]
$$
we see that (\ref{differ}) holds for $C_{12}(\lda,\mu)$ given in
(\ref{defC})\finprf

\begin{theo}
The commutation relations for the quantum matrices
$\ct^\pm(\lda;x)$ and $\ct(\lda)$ for real $\lda$ and $\mu$ take
the following form:
\begin{eqnarray}
\calr_{12}(\lda-\mu)\ct^-_1(\lda;x)\ct^-_2(\mu;x)C_{12}(\lda,\mu)
&=&\ct^-_2(\mu;x)\ct^-_1(\lda;x)C_{21}(\mu,\lda)\calr_{12}(\lda-\mu)
\label{RTT-}\nonu
\calr_{12}(\lda-\mu)C_{12}(\mu,\lda)\ct^+_1(\lda;x)\ct^+_2(\mu;x)
&=&C_{21}(\lda,\mu)\ct^+_2(\mu;x)\ct^+_1(\lda;x)\calr_{12}(\lda-\mu)
\label{RTT+}\nonu
\calr_{12}^+(\lda-\mu)\ct_1(\lda)\ct_2(\mu)
&=&\ct_1(\mu)\ct_2(\lda)\calr_{12}^-(\lda-\mu) \label{R+TTR-}
\end{eqnarray}
where, defining $\1_K=\sum_{i=1}^K E_{ii}$,
\begin{eqnarray*}
\calr_{12}^\pm(\lda-\mu)&=&\frac{-ig}{(\lda-\mu)}\1_K\otimes \1_K
+P_{12}+\pi_{21}\\
&+&\frac{(\lda-\mu)^2+g^2}{(\lda-\mu+i\varepsilon)^2}\pi_{12}
+\frac{\lda-\mu-ig}{\lda-\mu}E_{K+1,K+1}\otimes
E_{K+1,K+1}\\
&\pm&\pi g\delta(\lda-\mu)\Big(\1_K\otimes E_{K+1,K+1}-
 E_{K+1,K+1}\otimes \1_K\Big)
\end{eqnarray*}
\end{theo}

\prf We start with the proof of the first equality. Lemma
\ref{lemlambda} gives
\begin{equation*}
\lim_{y\to
-\infty}\ct_1(\lda;x,y)\,\ct_2(\mu;x,y)\,\xi_{12}(\lda,\mu;y)=
\ct^-_1(\lda;x)\,\ct^-_2(\mu;x)C_{12}(\lda,\mu)
\end{equation*}
which in turn yields
\begin{equation*}
\lim_{y\to-\infty}\ct_2(\mu;x,y)\ct_1(\lda;x,y)\xi_{21}(\mu,\lda;y)
=\ct^-_2(\mu;x)\ct^-_1(\lda;x)C_{21}(\mu,\lda)
\end{equation*}
Multiplying (\ref{finitevolRTT}) on the
right by $\xi_{12}(\lda,\mu;y)$ and using the property
$$
\calr_{12}(\lda-\mu)\xi_{12}(\lda,\mu;y)=\xi_{21}(\mu,\lda;y)\calr_{12}(\lda-\mu)
$$
we get
$$
\calr_{12}(\lda-\mu)\ct_1(\lda;x,y)\ct_2(\mu;x,y)\xi_{12}(\lda,\mu;y)
=\ct_2(\mu;x,y)\ct_1(\lda;x,y)\xi_{21}(\mu,\lda;y)\calr_{12}(\lda-\mu)
$$
which gives the first equality in the limit $y\to -\infty$. The second
equality
is proved along the same line of argument. Now, combining
the two equations and using the properties
\begin{equation*}
\ct(\lda)=\ct^+(\lda;x)\ct^-(\lda;x)~~\text{and}~~\ct^+_2(\mu;x)\ct^-_1(\lda;x)=\ct^-_1(\lda;x)\ct^+_2(\mu;x)
\end{equation*}
we get
\begin{eqnarray*}
\calr_{12}(\lda-\mu)C_{12}(\mu,\lda)\ct_1(\lda)\ct_2(\mu)C_{12}(\lda,\mu)
~~~~~~~~~~~~~~~~~~~~~~~~~~~~~~~~~~~~~~~\nonu
~~~~~~~~~~~~~~~~~~~~~~=C_{21}(\lda,\mu)\ct_2(\mu)\ct_1(\lda)C_{21}(\mu,\lda)\calr_{12}(\lda-\mu)
\end{eqnarray*}
which take the form (\ref{R+TTR-}) if we define
\begin{eqnarray}
\calr_{12}^+(\lda-\mu)&=& C_{12}^{-1}(\lda,\mu)\Pi_{12}\calr_{12}(\lda-\mu)C_{12}(\mu,\lda)\\
\calr_{12}^-(\lda-\mu)&=&C_{12}(\mu,\lda)\Pi_{12}\calr_{12}(\lda-\mu)C_{12}^{-1}(\lda,\mu)
\end{eqnarray}
whose explicit calculation we leave to the reader.
\finprf

Let us
extract the information contained in (\ref{R+TTR-}). We start by
particularizing some entries of the quantum monodromy matrix
($i,j=1,\ldots,K$):
\begin{eqnarray}
t_{ij}(\lda)&=&(\ct(\lda))_{ij}\\
b_j(\lda)&=&(\ct(\lda))_{j,K+1}\\
D(\lda)&=&(\ct(\lda))_{K+1,K+1}
\end{eqnarray}

\begin{theo}\label{exrel}
The exchange relations of the entries of the quantum
monodromy matrix read as follows
\begin{eqnarray}
\left[\hspace{-4pt}\left[\rule{0ex}{2.4ex}t_{ij}(\lda),t_{kl}(\mu)\right]\hspace{-4pt}\right]
&=&i
{g}\,(-1)^{[j][k]+[i]([j]+[k])}\,\frac{t_{kj}(\lda)t_{il}(\mu)
    -t_{kj}(\mu)t_{il}(\lda)}{\lda-\mu}\label{quantumSY}\\
t_{ij}(\lda)D(\mu)&=&D(\mu)t_{ij}(\lda)\label{symSY}\\
D(\lda)D(\mu)&=&D(\mu)D(\lda)\label{involution}\\
b_j(\lda)b_k(\mu)&=&\frac{\mu-\lda}{\mu-\lda-ig}
(-1)^{jk}b_k(\mu)b_j(\lda)-\frac{ig}{\mu-\lda-ig}b_j(\mu)b_k(\lda)
\qquad\label{bb}\\
b_j(\lda)D(\mu)&=&\frac{\lda-\mu-ig}{\lda-\mu-i\varepsilon}D(\mu)b_j(\lda)\label{bd}
\end{eqnarray}
\end{theo}

\prf By direct calculation\finprf

Relations (\ref{quantumSY}-\ref{involution}) are the quantum
counterparts of eqs (\ref{SY}-\ref{DD}) and the same
interpretation holds but for the quantum hierarchy here. As such,
the super-Yangian $Y(gl(M|N))$ is a quantum symmetry of the
hierarchy generated by $D(\lda)$, which is just the quantum analog
of Property \ref{firsthami} as can be seen from
\begin{eqnarray*}
D(\lda)&=&1+\frac{ig}{\lda}N+\frac{g}{\lda^2}\left(P-\frac{g}{2}N(N-1)\right)\\
&
&+\frac{ig}{\lda^3}\left(H+g(N-1)P-\frac{g^2}{6}N(N-1)(N-2)\right)+O\left(\frac{1}{\lda^4}\right)
\end{eqnarray*}

\subsection{ZF algebra from Lax pair}
The two relations (\ref{bb})and (\ref{bd}) will allow us to
recover the ZF algebra. Indeed, all the quantities of Theorem
\ref{exrel} are functionals of $\Phi, \Phidag$, themselves
involving the ZF generators (cf (\ref{quantum-field-comp})), and
one can get the ZF algebra out of them as follows
\begin{prop}\label{prop-a*a*}
Defining $a_j(\lda)=\frac{1}{\sqrt{\pi g}}b_j(\lda)D(\lda)^{-1}$,
eqs (\ref{bb}) and (\ref{bd}) give
\begin{equation}
a_j(\lda)a_k(\mu)=\frac{\mu-\lda}{\mu-\lda+ig}
(-1)^{jk}a_k(\mu)a_j(\lda)
-\frac{ig}{\mu-\lda+ig}a_j(\mu)a_k(\lda)
\end{equation}
\end{prop}
\prf Direct calculation from Theorem \ref{exrel}.\finprf

To
complete our algebra, we need the exchange relations between
$a_j(\lda)$ and $\adag_k(\mu)$. Contrary to the original one
(bosonic) component case, this is not directly obtained from what
we already have since there is no simple conjugate relationship
for the entries of the monodromy matrix. We are naturally led to
introduce a conjugate Lax super-matrix defined by:
\begin{equation}
\bL(\lda;x)=-\frac{i\lda}{2}\Sigma-i\sqrt{g}\phidag_j(x)
E_{K+1,j}+i\sqrt{g}\phi_j(x) E_{j,K+1}
\end{equation}
and the associated transition matrix
\begin{equation}
\prt_x\bct(\lda;x,y)=\ddagger\bct(\lda;x,y)\bL(\lda;x)\ddagger
\end{equation}
Now, to obtain information between the entries of $\ct(\lda;x,y)$
and $\bct(\mu;x,y)$ following the same steps as in lemmas
\ref{deriveTT}-\ref{RL} and Theorem \ref{finitevolRTTtheo}, one
sees that we actually need to work with the super-transposed Lax
matrix. The corresponding operation on an even super-matrix
$A=\sum_{i,j=1}^{K+1}A_{ij}E_{ij}$ reads
\begin{equation}
A^t=\sum_{i,j=1}^{K+1}A_{ij}\,E_{ij}^t=\sum_{i,j=1}^{K+1}(-1)^{[i]([i]+[j])}A_{ji}E_{ij}
\end{equation}
It satisfies $(A^t)^t=A$ and $(AB)^t=B^tA^t$ for any even
super-matrices $A$ and $B$. We get
\begin{equation}
L^t(\lda;x)=\frac{i\lda}{2}\Sigma+i\sqrt{g}(-1)^{[j]}\phi_j(x)
E_{K+1,j}-i\sqrt{g}\phidag_j(x) E_{j,K+1}
\end{equation}
and the associated transition matrix
\begin{equation}
\prt_x\ct^t(\lda;x,y)=\ddagger\ct^t(\lda;x,y) L^t(\lda;x) \ddagger
\end{equation}
Therefore, instead of (\ref{deriveTT1}) we get
\begin{eqnarray}
\prt_x \{\bct_1(\lda;x,y)\ct_2^t(\mu;x,y)\}&=&\ddagger
\bct_1(\lda;x,y)\ct^t_2(\mu;x,y)\Gamma_{12}(\lda,\mu;x)\ddagger\\
\prt_x \{\ct_1^t(\mu;x,y)\bct_2(\lda;x,y)\}&=&\ddagger
\ct_1^t(\mu;x,y)\bct_2(\lda;x,y)\Gamma'_{12}(\lda,\mu;x)\ddagger
\end{eqnarray}
with
\begin{eqnarray*}
\Gamma_{12}(\lda,\mu;x)&=&\bL_1(\lda;x)+L^t_2(\mu;x)+g\pi_{12}^{t_2}\\
\Gamma'_{12}(\lda,\mu;x)&=&L^t_1(\mu;x)+\bL_2(\lda;x)+g\pi_{12}^{t_1}
\end{eqnarray*}
Now the key point is to find an invertible numerical matrix
$\calr'_{12}(\lda)$ solution of the new Yang-Baxter equation
$$\calr'_{12}(\lda,\mu)\Gamma_{12}(\lda,\mu;x)=\Gamma_{21}(\lda,\mu;x)\calr'_{12}(\lda,\mu)$$
It is given by
\begin{equation}
\calr'_{12}(\lda,\mu)=\frac{ig}{\lda-\mu}\Pi_{12}^{t_1}+
\frac{\lda-\mu-ig(M-N)}{\lda-\mu}\Pi_{12}
\end{equation}
Following the same procedure as above, we finally deduce the
infinite volume commutation relations under the form
\begin{equation}
\label{R+TTR-dag}
{\calr'}^+_{12}(\lda-\mu)\bct_1(\lda)\ct^t_2(\mu)=
\ct_1^t(\mu)\bct_2(\lda){\calr'}^-_{12}(\lda-\mu)
\end{equation}
with
\begin{eqnarray*}
{\calr'}^\pm_{12}(\lda-\mu)&=&\frac{ig}{\lda-\mu}P_{12}^{t_1}
+\frac{\lda-\mu-ig(M-N)}{\lda-\mu}\left(P_{12}+\pi_{12}+\pi_{21}\right)\\
&+&\frac{(\lda-\mu-ig)(\lda-\mu-ig(M-N))}{(\lda-\mu+i\varepsilon)^2}E_{K+1,K+1}\otimes
E_{K+1,K+1} \\
 &\mp&\pi g\delta(\lda-\mu)\Big(\pi_{21}^{t_1}-\pi_{21}^{t_2}\Big)
\end{eqnarray*}
All these results are the generalization to the graded case of
\cite{PWZ} ($K$, the total number of bosonic or fermionic
particles is replaced in our case by $M-N$, the difference of
bosonic and fermionic particles). Accordingly, we get the same
conclusions collected in the following proposition
\begin{prop}
Let $\adag_i(\lda)=\frac{1}{\sqrt{\pi
g}}(D^{-1})^\dagger(\lda)b^\dagger_j(\lda)$, then:
\begin{eqnarray}
a_i(\lda)\adag_j(\mu)&=&
\frac{\lda-\mu}{\lda-\mu+ig}(-1)^{[i][j]}\adag_j(\mu)a_i(\lda)\nonu
&
&-\delta_{ij}\frac{ig}{\lda-\mu+ig}\sum_{\ell=1}^K\adag_{\ell}(\mu)a_{\ell}(\lda)
+\delta_{ij}\delta(\lda-\mu)
\label{aadag}\\
\adag_i(\lda)\adag_j(\mu)
&=&\frac{\mu-\lda}{\mu-\lda+ig}(-1)^{[i][j]}\adag_j(\mu)\adag_i(\lda)
-\frac{ig}{\mu-\lda+ig}\adag_i(\mu)\adag_j(\lda) \label{aa}
\end{eqnarray}
\end{prop}

\prf Noting that
$$b_j(\lda)=\ct^t(\lda)_{K+1,j},~~
D(\lda)=\ct^t(\lda)_{K+1,K+1}$$
$$b^\dagger_j(\lda)=\ct^\dagger(\lda)_{K+1,j},~~
D^\dagger(\lda)=\ct^\dagger(\lda)_{K+1,K+1}$$ (\ref{R+TTR-dag})
gives
\begin{eqnarray*}
D^\dagger(\lda)D(\mu)&=&D(\mu)D^\dagger(\lda)\\
D(\mu)b^\dagger_i(\lda)&=&\frac{\lda-\mu-ig}{\lda-\mu+i\varepsilon}b^\dagger_i(\lda)D(\mu),~~
b_j(\mu)D^\dagger(\lda)=\frac{\lda-\mu-ig}{\lda-\mu+i\varepsilon}D^\dagger(\lda)b_j(\mu)\\
b_i(\lda)b^\dagger_j(\mu)&=&\frac{\mu-\lda-ig}{\mu-\lda+i\varepsilon}(-1)^{[i][j]}b^\dagger_j(\mu)b_i(\lda)\nonu
&
&+\delta_{ij}\frac{ig(\mu-\lda-ig)}{(\mu-\lda+i\varepsilon)^2}\sum_{\ell=1}^K
b^\dagger_{\ell}(\mu)b_{\ell}(\lda)+\delta_{ij}\pi
g\delta(\lda-\mu)D^\dagger(\mu)D(\lda)
\end{eqnarray*}
which in turn yields (\ref{aadag}). The proof of (\ref{aa}) is
similar. \finprf

\section{Explicit construction of the super-Yangian generators\label{sec4}}
\subsection{Super-Yangian generators in terms of canonical
fields\label{canodescription}}
 We consider the classical case. The
quantum case can be done in a similar way, with correction terms
due to the non-commutativity of the fields $\Phi$, $\Phidag$.

For any $K\times K$-matrix $\sigma\in gl(M|N)$, we introduce:
\begin{eqnarray}
    Q_{\sigma}^{(0)} &=& \int dx\, \Phidag(x)\,\sigma\,\Phi(x) \ =\ \int
    dx\, \sum_{j,k=1}^{K}\phidag_{j}(x)\,\sigma^{jk}\,\phi_{k}(x)
    \label{Q0sig}\\
    Q_{\sigma}^{(1)} &=& \int dx\,
    \Phidag(x)\,\sigma\,\prt\Phi(x) -\frac{g}{2} \int
    dxdy\, sg(x-y)\,\Phidag(x)\sigma\Phi(y) \,\cdot\, \Phidag(y)\Phi(x)
    \label{Q1sig}\\
    Q_{\sigma}^{(2)} &=& \int dx\,
    \Phidag(x)\,\sigma\,\prt^2\Phi(x) \label{Q2sig}\\
    && -\frac{g}{2}\int
    dxdy\, sg(x-y)\,\left(\rule{0ex}{2.4ex}\Phidag(x)\sigma\prt\Phi(y)-
    \prt\Phidag(x)\sigma\Phi(y)\right)\,\Phidag(y)\Phi(x)\nonu
    && +\frac{g^2}{4}\int  dxdydz\, sg(x-y)sg(y-z)\;\Phidag(y)\Phi(x)\,\cdot\,
    \Phidag(x)\sigma\Phi(z)\,\cdot\, \Phidag(z)\Phi(y)\nonumber
\end{eqnarray}
The coefficients in (\ref{Q1sig}) and (\ref{Q2sig}) are fixed in such
a way that
\begin{equation}
\{H,Q_{\sigma}^{(n)}\} = 0,\ n=0,1,2
\end{equation}
so that $Q_{\sigma}^{(n)}$ are indeed symmetry generators of the NLSS
eq.
With these definitions, it is a simple calculation to show:
\begin{eqnarray}
    \{Q_{\sigma}^{(0)},Q_{\omega}^{(n)}\} &=& i
    Q_{\lcom\sigma,\omega\rcom}^{(n)}\qquad n=0,1,2
    \label{Q0Qn}\\
        \{Q_{\sigma}^{(1)},Q_{\omega}^{(1)}\} &=& i
    Q_{\lcom\sigma,\omega\rcom}^{(2)}-i \left(-\frac{g}{2}\right)^2\int dx\,dy\,dt\,
    S(x,y,t)\Big(
    \Phidag(x)\sigma\Phi(y)\cdot\Phidag(y)\omega\Phi(t)\nonu
    &&\qquad\qquad-
    \Phidag(x)\omega\Phi(y)\cdot\Phidag(y)\sigma\Phi(t)\Big)
    \Phidag(t)\Phi(x)
    \label{Q1Q1}\\
    S(x,y,t) &=& sg(t-x)sg(x-y) + sg(x-y)sg(y-t) + sg(y-t)sg(t-x)
    \nonumber
\end{eqnarray}
Equation (\ref{Q0Qn}) shows that $Q_{\sigma}^{(0)}$, $\sigma\in gl(M|N)$, generates a
$gl(M|N)$ superalgebra, and that $Q_{\sigma}^{(n)}$ ($n$ fixed) form a
representation of it. The second term in (\ref{Q1Q1}) reflects the
non-linear commutation relation of the super-Yangian.

Note that we have
\begin{equation}
Q^{(0)}_{{\mathbb I}} = N \mb{and} Q^{(1)}_{{\mathbb I}} = P
\end{equation}
so that eq. (\ref{Q0Qn}) shows that
$Q_{\sigma}^{(n)}$ commutes with $N$ and $P$. Moreover, we
have the supersymmetry-like relations:
\begin{equation}
\begin{array}{l}
\{Q_{\sigma}^{(0)},Q_{\sigma}^{(0)}\}=2iN\\
\{Q_{\sigma}^{(0)},Q_{\sigma}^{(1)}\}=2iP
\end{array} \mb{as soon as} \sigma^2=\II
\mbox{ and }[\sigma]=1
\end{equation}
However, let us remark that $Q^{(2)}_{{\mathbb I}}$ is not the NLSS
Hamiltonian:
\begin{equation*}
Q^{(2)}_{{\mathbb I}} = H +
\frac{g^2}{4}\int  dxdydz\, sg(x-y)sg(y-z)\;\Phidag(y)\Phi(x)\,\cdot\,
    \Phidag(x)\sigma\Phi(z)\,\cdot\, \Phidag(z)\Phi(y)
\end{equation*}
$Q^{(2)}_{{\mathbb I}}$ corresponds to a central generator which,
if it were the Hamiltonian, would lead to non-local equation
of motion for $\Phi$.
On the contrary, $H$ commutes with the generators $Q^{(n)}_{\sigma}$
and provides local equation of motion.
\subsection{Super-Yangian generators in terms of ZF generators}

We have obtained the ZF-algebra
(\ref{ZF-algebra1},\ref{ZF-algebra3}) from the commutation
relations of the quantum monodromy matrix. This shows the
central importance of this algebra and one is naturally led to
take it as a starting point. This is the very idea developed in
\cite{Rag} and we use it to construct a realization of the
generators of the super-Yangian symmetry in terms of the ZF
oscillators.

First of all, we need to generalize all the basic results of
\cite{Rag} to our graded formalism. It is actually readily
obtained since the fundamental idea of the properties given in
\cite{Rag,MRSZ} is the possibility of relabelling the auxiliary
spaces
 which holds for our global
formalism as the reader can check. Thus, we are in position to
apply any result from \cite{Rag} in our context. Here is our
strategy : we start from the ZF algebra (corresponding to
the algebra $\ca_R$ in \cite{Rag}), introduce the associated
well-bred vertex operator $T(\lda)$ and use the explicit
expression of our $R$-matrix to derive the first two
terms of the expansion of $T(\lda)$ in power series of
$\lda^{-1}$. Then we show that this approach actually coincides
with the previous Lax pair formulation so that we have a
realization of the generators of the super-Yangian symmetry for
the hierarchy associated to the nonlinear super-Schr\"odinger
equation in terms of the ZF oscillators. This completes and
confirms the deep relationships between the quantum canonical
field description (cf section \ref{canodescription}) and the ZF
algebra approach.

\begin{defi}
The vertex operators $T_{ij}(\lda),~~(i,j=1,\ldots,K)$
associated to the ZF algebra $\ca_R$ are defined by
$T(\lda)=T^{ij}(\lda)E_{ij}\in\ca_R\otimes \CC^{K^2}$ with
\begin{equation}
T_\infty(\lda)=\1 +\sum_{n=1}^\infty
\frac{(-1)^{n+1}}{n!}\adag_{n\ldots1}T^{(n)}_{\infty 1\ldots
n}a_{1\ldots n}
\end{equation}
where
\begin{eqnarray*}
\adag_{n\ldots1}&=&(a_{1\ldots
n})^\dagger=\adag_n(k_n)\ldots\adag_1(k_1)\\
T^{(n)}_{\infty 1\ldots n}&=&T^{(n)}_{\infty 1\ldots
n}(\lda,k_1,\ldots,k_n)\in(\CC^{\otimes
K^2})^{\otimes(n+1)}(\lda,k_1,\ldots,k_n)
\end{eqnarray*}
and integration is implied over the spectral parameters
$k_1,\ldots,k_n$ (the summation over the auxiliary spaces being
understood as in the appendix).

$T_\infty(\lda)$ is said to be well-bred (on $\ca_R$) if
\begin{equation}
\label{bred}
T_\infty(\lda)a_1(\mu)=R_{1\infty}(\mu-\lda)a_1(\mu)T_\infty(\lda)~~and~~
T_\infty(\lda)\adag_1(\mu)=\adag_1(\mu)R_{\infty
1}(\lda-\mu)T_\infty(\lda)\quad
\end{equation}
with $R$ given by (\ref{R-matrix}).
\end{defi}
Then, from \cite{Rag} we can directly assert
\begin{prop}
The well-bred vertex operators $T_\infty(\lda)$ obey
Faddeev-Reshetikhin-Takhtajan (FRT) relations
\begin{equation}
\label{FRT}
R_{\infty\infty'}(\lda-\mu)T_\infty(\lda)T_{\infty'}(\mu)=T_{\infty'}(\mu)T_\infty(\lda)R_{\infty\infty'}(\lda-\mu)
\end{equation}
so that they generate the super-Yangian algebra $Y(gl(M|N))$. In
addition, they form a symmetry super-algebra for the hierarchy
$H^{(n)}$ defined by
\begin{equation}
H^{(n)}=\int_{-\infty}^\infty
dk~k^n\adag(k)a(k)~~~~~~~~n\in\ZZ_+
\end{equation}
forming an Abelian algebra of hermitian operators and governing
the flows of the scattering operators $a,\adag$ as follows:
\begin{eqnarray*}
e^{iH^{(n)}t}a(k)e^{-iH^{(n)}t}=e^{-ik^{n}t}a(k)\\
e^{iH^{(n)}t}\adag(k)e^{-iH^{(n)}t}=e^{ik^{n}t}\adag(k)
\end{eqnarray*}
\end{prop}

Now, recalling the results obtained in Section \ref{timeevo},
Property \ref{firsthami} and
eqs.(\ref{quantumSY}-\ref{involution}), we see that both
descriptions of our integrable system (in terms of canonical
fields or ZF scattering operators) are equivalent. But in this
operation, we have gained an explicit realization of the
super-Yangian generators.

To do this, we use the inductive relations obtained in Theorem 3.3
of \cite{Rag} order by order in the spectral parameter $\lda$. Let
us rewrite
\begin{equation}
T_\infty(\lda)=\1 +\frac{ig}{\lda}\sum_{p=0}^\infty
T_\infty^{\{p\}}\lda^{-p}
\end{equation}
where accordingly,
\begin{equation*}
T_\infty^{\{p\}}=\sum_{n=1}^\infty
\frac{(-1)^{n+1}}{n!}\adag_{n\ldots1}T^{(n)\{p\}}_{\infty 1\ldots
n}a_{1\ldots n}
\end{equation*}
for some $T^{(n)\{p\}}_{\infty 1\ldots n}\in(\CC^{\otimes
K^2})^{\otimes(n)}(k_1,\ldots,k_n)$.

Our goal is to determine $T_\infty^{\{0\}}$ and
$T_\infty^{\{1\}}$, that is the first two "levels" of the
super-Yangian generators. To do this we note that the inductive
relations of Theorem 3.3 in \cite{Rag} at first order in $\lda$
take the form
\begin{equation}
T^{(n+1)}_{\infty 0\ldots n}=T^{(n)}_{\infty 1\ldots
n}-T^{(n)}_{\infty 0\ldots n-1}+O(\lda^{-2})
\end{equation}
which, under the knowledge of
$$
T^{(1)\{0\}}_{\infty 0}=\1 +P_{\infty0}
$$ 
yields
$$
T^{(n+1)\{0\}}_{\infty 0\ldots n}=(-1)^n  \sum_{k=0}^n(-1)^k
\bino{n}{k}P_{\infty k}
$$ 
where $P_{ij}$ is the super-permutation 
of auxiliary spaces $i$ and $j$, so that
\begin{equation}
\label{T0}
T_\infty^{\{0\}}=\sum_{n=0}^\infty\frac{(-1)^{n+1}}{n!}\sum_{k=0}^n(-1)^{n-k}\bino{n}{k}\adag_{n\ldots
0}P_{\infty k}a_{0\ldots n}
\end{equation}
Now that we have the explicit form of $T_\infty^{\{0\}}$ we can use
it to evaluate the commutator $[T_{\infty'}^{\{1\}},T_\infty^{\{0\}}]$
directly and compare the result to that obtained from the FRT
relations (\ref{FRT}) at order $\lda^{-2}$. The latter calculation
yields
\begin{equation}
\label{T1T0FRT}
[T_{\infty'}^{\{1\}},T_\infty^{\{0\}}]=[P_{\infty'\infty},T_\infty^{\{1\}}]
\end{equation}
As for the former, the well-bred relations (\ref{bred}) at order
$\lda^{-2}$ read
\begin{eqnarray*}
\left[T_\infty^{\{0\}},a_0(\mu)\right]&=&(\1 +P_{0\infty})a_0(\mu)\\
\left[T_\infty^{\{1\}},a_0(\mu)\right]&=&\mu(\1
+P_{0\infty})a_0(\mu)+ig(\1 +P_{0\infty})a_0(\mu)(\1
+T_{\infty}^{\{0\}})\\
\left[T_\infty^{\{0\}},\adag_0(\mu)\right]&=&-\adag_0(\mu)(\1 +P_{\infty 0})\\
\left[T_\infty^{\{1\}},\adag_0(\mu)\right]&=&-\mu\adag_0(\mu)(\1
+P_{\infty 0})+ig\adag_0(\mu)(\1 +P_{\infty 0})(\1
-T_{\infty}^{\{0\}})
\end{eqnarray*}
which will be useful in calculating
$$[T_{\infty'}^{\{1\}},T_\infty^{\{0\}}]=\sum_{n=0}^\infty\frac{(-1)^{n+1}}{n!}
\sum_{k=0}^n(-1)^{n-k}\bino{n}{k}[T_{\infty'}^{\{1\}},\adag_{n\ldots
0}P_{\infty k}a_{0\ldots n}]$$

Note that this procedure can  be iterated to evaluate
$T_{\infty'}^{\{n\}}$ for an arbitrary $n$ through
$[T_{\infty'}^{\{n\}},T_\infty^{\{0\}}]$. Now,
\begin{eqnarray*}
[T_{\infty'}^{\{1\}},\adag_{n\ldots 0}P_{\infty k}a_{0\ldots
n}]&=&\sum_{i=0}^n\adag_{n\ldots 0}P_{\infty
k}a_0\ldots [T_{\infty'}^{\{1\}},a_i]\ldots a_n
\\
& &+\sum_{i=0}^n\adag_n\ldots
[T_{\infty'}^{\{1\}},\adag_i]\ldots\adag_0 P_{\infty k}a_{0\ldots
n}\\
&=&[P_{\infty\infty'},
(\mu_k-n-1)\adag_{n\ldots 0}P_{\infty k}a_{0\ldots n}]
\\
& &+\adag_{n\ldots 0}[P_{\infty\infty'},P_{\infty
k}]T_{\infty'}^{\{0\}}a_{0\ldots n}
+\sum_{i=0}^{k-1}\adag_{n\ldots 0}[P_{\infty
k},P_{\infty\infty'}]P_{\infty'i}a_{0\ldots n}\\
& &+\sum_{i=k+1}^n\adag_{n\ldots 0}P_{\infty'i}[P_{\infty
k},P_{\infty\infty'}]a_{0\ldots n}
\end{eqnarray*}

This expression can be considerably simplified in
$[T_{\infty'}^{\{1\}},T_\infty^{\{0\}}]$ using the properties of the
binomial coefficients to combine the last three terms. Inserting
(\ref{T0}) and using the property
\begin{equation*}
\sum_{n=k}^{i-1}\bino{N}{n}\alpha_k^n\alpha_{i-n}^{N-n}=\alpha_k^N-\alpha_i^N
\mb{where} \alpha_k^n=(-1)^{k-1}\bino{n-1}{k-1}
\end{equation*}
proved in \cite{MRSZ}, we get (after a convenient relabelling of
the auxiliary spaces):
\begin{equation*}
[T_{\infty'}^{\{1\}},T_\infty^{\{0\}}]=
\Big[P_{\infty\infty'},\sum_{n=1}^\infty
\frac{(-1)^n}{n!}\sum_{k=1}^n\alpha_k^n\adag_{1\ldots
n}\lbrace(\mu_k-ign)P_{\infty k}-ig\sum_{i=1}^{k-1}P_{\infty
i}P_{\infty k}\rbrace a_{n\ldots 1}\Big]
\end{equation*}
Comparing this last expression with (\ref{T1T0FRT}), we get the
explicit form for $T_{\infty}^{\{1\}}$ (up to a term
proportional to $\II_\infty$).

To conclude, we can recast this expression as
\begin{eqnarray}
T_{\infty}^{\{1\}}&=&\sum_{n=1}^\infty\frac{(-1)^n}{n!}\sum_{k=1}^n\alpha_k^n\,
\adag_{1\ldots n}\Big(\mu_k P_{\infty k} -ig
\sum_{i=1}^{k-1} P_{\infty k}P_{\infty i}\Big)
a_{n\ldots1}+ig\,T_{\infty}^{\{0\}}T_{\infty}^{\{0\}}\qquad\quad
\end{eqnarray}
In the case of $gl(N)$, we recover the results of \cite{MRSZ},
although in a different basis:
\begin{eqnarray}
T_{ij}^{\{0\}} &=& \sum_{n=0}^\infty\frac{(-1)^{n+1}}{n!}
\sum_{k=0}^n\alpha_k^n\,\adag_{n\ldots0}\,E_{ji}^{(k)}\,a_{0\ldots n}
\nonu
T_{ij}^{\{1\}}&=&\sum_{n=1}^\infty\frac{(-1)^n}{n!}\sum_{k=1}^n\alpha_k^n\,
\adag_{1\ldots n}\Big(\mu_k E_{ji}^{(k)} -ig
\sum_{\ell=1}^{k-1} \sum_{m=1}^{N}E^{(\ell)}_{jm}E^{(k)}_{mi}\Big)
a_{n\ldots1}+ig\Big(T^{\{0\}}\Big)^{2}_{ji}
\nonumber
\end{eqnarray}
where $E^{(\ell)}_{ij}$ denotes the $E_{ij}$ matrix in the $\ell^{th}$
auxiliary space.

For $gl(M|N)$, similar formulae may also be obtained, taking care of the
$\ZZ_{2}$-graded tensor products.

\section*{Conclusion}
We solved a vectorial version of the Nonlinear Schr\"odinger
equation which contains fermions and bosons at the same time. 
We first introduced it classically using a
$\ZZ_2$-graded formalism. At the quantum level, special attention
was paid to the resolution using a super ZF algebra associated to
the $R$-matrix of the super-Yangian $Y(gl(M|N))$. The
integrability and symmetry of our system was studied through a Lax
pair formalism and it is worth stressing the deep interplay
between canonical and (ZF) algebraic formalisms. The ZF algebra 
allowed us to compute the correlation functions. 
Further investigations can
be performed in this direction to study super-versions of known
integrable systems. One can also study these super-versions when
 a boundary is introduced, using generalizations of the ZF
algebra (boundary algebras).

\appendix
\section{Appendix}
\subsection{Auxiliary spaces}
We define in the auxiliary space, a $K$-column vector
 $e_{j}$  with 1 at row $j$ and 0
elsewhere, its transpose, the row vector  $\edag_i=(0,\ldots,1,\ldots,0)$
and the matrices $E_{ij}$, with 1 at position $(i,j)$.

Here and below, the vectors $e_{i}$, $e^\dagger_{i}$, and the
matrices $E_{ij}$  will be $\ZZ_{2}$-graded:
$$
[e_{i}]=[e^\dagger_{i}]=[i]\ ;\ [E_{ij}]=[i]+[j] \mb{with} [i]=
\begin{cases}
0 &\text{for}~i=1,...,M \\
1 &\text{for}~i=M+1,...,N
\end{cases}
$$
Accordingly, the tensor product of auxiliary spaces will be also
$\ZZ_{2}$-graded, e.g.
$$
(\II\otimes e_{i})(E_{jk}\otimes \II)=(-1)^{[i]([j]+[k])}\,
E_{jk}\otimes e_{i}
$$
We will consider even objects in the following sense:
$v=v_{i}e_{i}$ and $U=U_{ij}E_{ij}$ (summation on repeated indices is understood)
are even iff $[v_{i}]=[i]$ and
$[U_{ij}]=[i]+[j]$. For example, the field $\Phi(x)$
 is even.

Note that, when dealing with tensor product of
auxiliary spaces,
 one has to be careful not to confuse (even) objects like
 $\gld_{1}=\gld\otimes \II=\sum_{i=1}^{K}\lda_{i}e_{i}\otimes \II$
 with their ($\ZZ_{2}$-graded) components $\lda_{i}$, $i=1,\ldots,K$.
 As a (tentative) clarifying notation, we will use boldface letters
 for the even objects, and ordinary letters for their components.

 Finally, in order to apply our formalism to derive the classical
NLSS equation, we will use  the  global Kronecker symbol:
\begin{equation}
\delta_{12}=\delta^{ij}(e_i\otimes
\edag_j)=(e_i\otimes\edag_i)
\end{equation}
and, accordingly
\begin{equation}
\delta_{21}=(-1)^{[i]}(\edag_i\otimes e_i)
\end{equation}

\subsection{Poisson Brackets}
For $F$ and $G$ two $(\Phi, \Phidag)$-functionals , their Poisson
bracket is defined by
\begin{equation}
\label{super-poisson-bracket}
\{F,G\}=i\sum_{\ell=1}^{K}\int_{-\infty}^{\infty}dx(-1)^{[F][\ell]}
\left((-1)^{[\ell]}\frac{\delta F}{\delta
\phi_\ell(x)} \frac{\delta G}{\delta\phidag_\ell(x)} -\frac{\delta
F}{\delta \phidag_\ell(x)} \frac{\delta
G}{\delta\phi_\ell(x)}\right)
\end{equation}
This bracket is a
graded Poisson bracket \ie it is bilinear, graded
antisymmetric, obeys the graded Leibniz rule and graded Jacobi identity.

To any graded PB, one can associate  a ``global'' Poisson
bracket,
defined for
the even functionals $\gF$ and $\gG$.
We introduce the notation $u_{\alpha}$ to denote either $e_{i}$
($\alpha=(0,i)$ and $[\alpha]=[i]$), $e^\dagger_{i}$
($\alpha=(i,0)$ and $[\alpha]=[i]$), or $E_{ij}$ ($\alpha=(i,j)$
and $[\alpha]=[i]+[j]$), so that any even object $\gF$ can be
written  $\gF=\sum_{\alpha}F_{\alpha}u_{\alpha}$ with
$[F_{\alpha}]=[\alpha]$.

On any even object, one defines the global PB
\begin{equation}
    \{\gF_{1},\gG_{2}\} = \sum_{\alpha,\beta}\{\gF_{\alpha},\gG_{\beta}\}
    \,u_{\alpha}\otimes u_{\beta} \label{globPB}
\end{equation}
It is bilinear, antisymmetric, obeys Leibniz rule and Jacobi identity.
Let us stress that this global PB is not graded (because of the use of
auxiliary spaces), but its "component" version indeed is graded.
\begin{lem}
    The  global PB (\ref{globPB}) corresponding to the graded PB
    (\ref{super-poisson-bracket}) can be rewritten as
\begin{equation}
    \{\gF_{1},\gG_{2}\} = i\int_{\RR}dx \left(
    \frac{\delta\,\gF_{1}}{\delta\Phi_{3/2}(x)}\,
\frac{\delta\,\gG_{2}}{\delta\Phidag_{3/2}(x)}\ -\
\frac{\delta\,\gG_{2}}{\delta\Phi_{3/2}(x)}\,
\frac{\delta\,\gF_{1}}{\delta\Phidag_{3/2}(x)}\right) \label{PB-sg}
\end{equation}
where we have introduced a third auxiliary space (labeled $3/2$)
which is "inserted" between the space 1 and the space 2. We have
also defined
\begin{equation}
\frac{\delta}{\delta\Phi(x)}=
\sum_{j=1}^{K}e^\dagger_{j}\,\frac{\delta}{\delta\phi_{j}(x)}
\mb{and} \frac{\delta}{\delta\Phidag(x)}=
\sum_{j=1}^{K}(-1)^{[j]}\,e_{j}\,\frac{\delta}{\delta\phidag_{j}(x)}
\end{equation}
\end{lem}
\prf Direct calculation.
\finprf

%\newpage

\end{document}